\documentclass{article}

\usepackage{amsmath, amssymb, mathtools}
\usepackage{cite}
\usepackage{color,soul}  
\usepackage{flafter}  
\usepackage{graphicx}
\usepackage{fullpage}
\usepackage[mathscr]{eucal}
\usepackage[small]{caption}
\usepackage{subfigure}
\usepackage{enumitem}
\usepackage{bm} 

\newcommand{\inputtikz}[1]{\includegraphics[scale=1]{image/#1}}

\newtheorem{thm}{Theorem}[section]

\newtheorem{rem}[thm]{Remark}

\let\originalleft\left
\let\originalright\right
\renewcommand{\left}{\mathopen{}\mathclose\bgroup\originalleft}
\renewcommand{\right}{\aftergroup\egroup\originalright}

\newcommand{\abs}[1]{\left\vert#1\right\vert}

\newcommand{\cell}{\mathcal{C}}
\newcommand{\cI}{\mathcal{I}}
\newcommand{\cK}{\mathcal{K}}
\newcommand{\cM}{\mathcal{M}}
\newcommand{\cR}{\mathcal{R}}
\newcommand{\cS}{\mathcal{S}}
\DeclareMathOperator{\diverg}{div}
\newcommand{\collFreq}{\nu}

\newcommand{\collOp}{\mathcal{Q}}
\newcommand{\Q}{\collOp}
\newcommand{\collOpMin}{\collOp^-}
\newcommand{\collOpPlus}{\collOp^+}
\newcommand{\collOper}[1]{\collOp\left(#1\right)(\v)}

\newcommand{\collOperPlus}[1]{\collOpPlus\left(#1\right)(\v)}
\newcommand{\collOpLin}{\mathcal{L}}

\newcommand{\Dtime}{\mathrm{D}_t}
\newcommand{\Dxv}{\mathrm{D}_{\xGrid,\vGrid}}
\newcommand{\Dt}{\Delta t}
\newcommand{\dt}{\delta t}
\newcommand{\dx}{\Delta x}
\newcommand{\fepsi}{f^\epsi}
\newcommand{\fSDSys}{\mathbf{f}} 		
\newcommand{\gepsi}{g^{\epsi}}
\newcommand{\hepsi}{h^{\epsi}}
\newcommand{\heatflux}{\mathbf{q}}
\newcommand{\hydroLimit}{\epsi \to 0}
\newcommand{\Ma}{\mathit{Ma}} 

\def\Ball{ {\cal B}}
\def\supp{{\rm Supp}}
\def\QL{\collOp^{R}}
\def\f{\hat f}
\def\bb{\hat \beta}
\def\B{\hat B}
\DeclareMathOperator{\sinc}{Sinc}

\newcommand{\Maxwellian}[2]{\cM_{#1}\left(#2\right)} 
\newcommand{\MaxArtScal}[1]{\Maxwellian{#1}{\rho^\epsi}} 		
\newcommand{\MaxBGKOneD}{\Maxwellian{v}{\fepsi}} 		
\newcommand{\MaxBGKMultiD}{\Maxwellian{\v}{\fepsi}}		
\newcommand{\MaxBGKLinMultiD}{\Maxwellian{\text{\scriptsize{lin}},\v}{\fepsi}} 		
\newcommand{\MaxGlobMultiD}{\cM^{\rhoinf,\vMacroMultiDinf,\Tinf}_{\v}} 

\newcommand{\pot}{\Phi}
\newcommand{\R}{\mathbb{R}}			
\newcommand{\rhoinf}{\rho^{\infty}}
\renewcommand{\S}{\mathbb{S}}
\newcommand{\Sdt}{S_{\dt}}

\newcommand{\sett}[3]{\left(#1\right)_{#2=1}^{#3}}
\newcommand{\sigmavec}{\boldsymbol{\sigma}}
\newcommand{\sten}{\mathcal{S}}
\newcommand{\Tinf}{T^{\infty}}
\newcommand{\ths}{\theta_\sigma}

\renewcommand{\v}{\mathbf{v}}
\renewcommand{\k}{\mathbf{k}}
\renewcommand{\l}{\mathbf{l}}
\newcommand{\m}{\mathbf{m}}
\newcommand{\vGrid}{\boldsymbol{v}}
\newcommand{\vMacroOneD}{\bar{v}}

\newcommand{\vMacroMultiD}{\mathbf{\vMacroOneD}}
\newcommand{\vMacroMultiDinf}{\vMacroMultiD^{\infty}}

\newcommand{\vpec}{\mathbf{c}}
\newcommand{\vrel}{\v_r}
\newcommand{\vstarb}{\v_{*}}
\newcommand{\x}{\mathbf{x}} 		
\newcommand{\xGrid}{\boldsymbol{x}}
\newcommand{\epsi}{\varepsilon}		

\begin{document}

\title{Projective and telescopic projective integration for the nonlinear BGK and Boltzmann equations}

\author{
Ward Melis
	\thanks{Department of Computer Science, K.U. Leuven, Celestijnenlaan 200A, 3001 Leuven, Belgium ({\tt ward.melis@cs.kuleuven.be}).} \and
Thomas Rey
	\thanks{Laboratoire Paul Painlev\'e, Universit\'e de Lille,	Cit\'e Scientifique, 59655 Villeneuve d'Ascq, France({\tt thomas.rey@math.univ-lille1.fr})} \and
Giovanni Samaey
	\thanks{Department of Computer Science, K.U. Leuven, Celestijnenlaan 200A, 3001 Leuven, Belgium ({\tt giovanni.samaey@cs.kuleuven.be}).}}
\maketitle

\begin{abstract}	
We present high-order, fully explicit projective integration schemes for nonlinear collisional kinetic equations such as the BGK and Boltzmann equation. The methods first take a few small (inner) steps with a simple, explicit method (such as direct forward Euler) to damp out the stiff components of the solution. Then, the time derivative is estimated and used in an (outer) Runge-Kutta method of arbitrary order. The procedure can be recursively repeated on a hierarchy of projective levels to construct telescopic projective integration methods. Based on the spectrum of the linearized collision operator, we deduce that the computational cost of the method is essentially independent of the stiffness of the problem: with an appropriate choice of inner step size, the time step restriction on the outer time step, as well as the number of inner time steps, is independent of the stiffness of the (collisional) source term. In some cases, the number of levels in the telescopic hierarchy depends logarithmically on the stiffness. We illustrate the method with numerical results in one and two spatial dimensions.
\end{abstract}

\noindent \textbf{Keywords:} \texttt{Boltzmann equation, BGK equation, Projective Integration, spectral theory, fast spectral scheme.}



\section{Introduction \label{sec:introduction}}
Kinetic equations represent a gas as a set of particles undergoing instantaneous collisions interspersed with ballistic motion \cite{CIP}. Nowadays, these models appear in a variety of sciences and applications, such as astrophysics, aerospace and nuclear engineering, semiconductors, fusion processes in plasmas, as well as biology, finance and social sciences. The common structure of such equations consists in a combination of a linear transport term with one or more interaction terms, which together dictate the time evolution of the distribution of particles in the (six-dimensional) position-velocity phase space. 
From a numerical point of view, it is clear that this results in a real challenge, since the computational cost immediately becomes prohibitive for realistic problems \cite{DimarcoPareschi2014}. Aside from the curse of dimensionality, there are many other difficulties which are specific to kinetic equations. 
We recall two among the most important ones. The first is the computational cost related to the evaluation of the collision operator, which implies the computation of multidimensional integrals in each point of the physical space \cite{FiMoPa2006, PaRuSINUM2000}. 
The second challenge is represented by the presence of multiple time scales in the collision dynamics, leading to a very small mean free path, at least in parts of the spatial domain. Usually, computational problems exhibit multiple regimes in different regions in space. This requires the development of adapted numerical schemes to avoid the resolution of the stiff dynamics \cite{DimarcoPareschi2013, Jin_review, Jin1999, BLM2008, degondrev}. 

Historically, two different approaches are generally used to tackle kinetic equations numerically: deterministic methods, such as finite volume, semi-Lagrangian and spectral schemes \cite{DimarcoPareschi2014}, and probabilistic methods, such as Direct Simulation Monte Carlo (DSMC) schemes \cite{bird, Caflisch98}.
Both methodologies have strengths and weaknesses. Deterministic methods can normally reach high orders of accuracy. Nevertheless, stochastic methods are often faster, especially
for solving steady problems, but, typically, exhibit lower convergence rates and difficulties in describing non-stationary and slow motion flows. In this paper, we will consider deterministic methods, in which we evaluate the collision operator using a fast spectral method, in the spirit of \cite{MoPa:2006}. For a comprehensive overview of numerical schemes for collisional kinetic equations, such as equation \eqref{eqCollision}, we refer to \cite{DimarcoPareschi2014} and references therein.

In this paper, we are specifically interested in the time discretization of kinetic equations with stiffness arising from multiple time scales in the collision operator. The stiffness is usually characterized by the (small) mean free path $\epsi$, and becomes infinite when $\epsi$ tends to zero. In that limit, a limiting macroscopic equation emerges in terms of a few moments of the particle distribution (density, momentum, energy); the full particle distribution then relaxes infinitely quickly to a Maxwellian distribution defined by these low-order moments. There is currently a large research effort in the design of algorithms that are uniformly stable in $\epsi$ and approach a scheme for the limiting equation when $\epsi$ tends to 0; such schemes are called asymptotic-preserving in the sense of Jin \cite{Jin1999}. Again, we refer to the recent review \cite{DimarcoPareschi2014} for a clear survey on numerical methods for kinetic equations. Here, we briefly review some achievements using different strategies. In \cite{Jin1999,Jin2000a}, separating the distribution function $f$ into its odd and even parts in the velocity variable results in a coupled system of transport equations where the stiffness appears only in the source term, allowing to use a time-splitting technique with implicit treatment of the source term; see also related work in \cite{Jin1999,Klar1999,Klar1999a}. Implicit-explicit (IMEX) schemes are an extensively studied technique to tackle this kind of problems, see \cite{ascher1995,FilbetJin2010} and references therein. Recent results in this setting were obtained by Dimarco et al.~to deal with nonlinear collision kernels \cite{DimarcoPareschi2013}, and an extension to hyperbolic systems in a diffusive limit is given in \cite{Boscarino2013}. 
A different approach, based on well-balanced methods, was introduced by Gosse and Toscani \cite{Gosse2003,Gosse2004}, see also \cite{Buet2007}. 
When the collision operator allows for an explicit computation, an explicit scheme can be obtained subject to a classical diffusion CFL condition by splitting the particle distribution into its mean value and the first-order fluctuations in a Chapman-Enskog expansion form \cite{Godillon-Lafitte2005}. Also closure by moments, e.g. \cite{Coulombel2005}, can lead to reduced systems for which time-splitting provides new classes of schemes \cite{Carrillo2008}. 
Alternatively, a micro-macro decomposition based on a Chapman-Enskog expansion has been proposed \cite{Lemou2008}, leading to a system of transport equations that allows to design a semi-implicit scheme without time splitting. A non-local procedure based on the quadrature of kernels obtained through pseudo-differential calculus was proposed in \cite{Besse2010}. 

A robust and fully explicit method, which allows for time integration of (two-scale) stiff systems with arbitrary order of accuracy in time, is projective integration. 
Projective integration was proposed in \cite{Gear2003projective} for stiff systems of ordinary differential equations with a clear gap in their eigenvalue spectrum. In such stiff problems, the fast modes, corresponding to the Jacobian eigenvalues with large negative real parts, decay quickly, whereas the slow modes correspond to eigenvalues of smaller magnitude and are the solution components of practical interest. Projective integration allows a stable yet explicit integration of such problems by first taking a few small (inner) steps using a step size $\dt$ with a simple, explicit method, until the transients corresponding to the fast modes have died out, and subsequently projecting (extrapolating) the solution forward in time over a large (outer) time step of size ${\Dt > \dt}$. In \cite{Lafitte2012}, projective integration was analyzed for kinetic equations with a diffusive scaling. An arbitrary order version, based on Runge-Kutta methods, has been proposed recently in \cite{LafitteLejonSamaey2015}, where it was also analyzed for kinetic equations with an advection-diffusion limit. In \cite{LafitteMelisSamaey2017}, the scheme was used to construct a explicit, flexible, arbitrary order method for general nonlinear hyperbolic conservation laws, based on relaxation to a kinetic equation. Alternative approaches to obtain a higher-order projective integration scheme have been proposed in \cite{Lee2007,Rico-Martinez}. 
These methods fit within recent research efforts on numerical methods for multiscale simulation \cite{E2007,Kevrekidis2003}.  

For problems exhibiting more than a single fast time scale, telescopic projective integration (TPI) was proposed \cite{Gear2003telescopic}. In these methods, the projective integration idea is applied recursively. Starting from an inner integrator at the fastest time scale, a projective integration method is constructed with a time step that corresponds to the second-fastest time scale. This projective integration method is then considered as the inner integrator of a projective integration method on yet a coarser level. By repeating this idea, TPI methods construct a hierarchy of projective levels in which each outer integrator step on a certain level serves as an inner integrator step one level higher. The idea was studied and tested for linear kinetic equations in \cite{MelisSamaey2017}. These methods turn out to have a computational cost that is essentially independent of the stiffness of the collision operator.  

We do not call projective integration methods asymptotic-preserving as such, because we cannot explicitly evaluate the scheme for $\epsi=0$ to obtain a classical numerical scheme for the limiting equation. Nevertheless, projective and telescopic projective integration methods share important features with asymptotic-preserving methods. In particular, their computational cost does (in many cases) not depend on the stiffness of the problem. To be specific, it was shown in \cite{MelisSamaey2017}, for linear kinetic equations, that the number of inner time steps at each level of the telescopic hierarchy is independent of the small-scale parameter $\epsi$, as is the step size of the outermost integrator. The only parameter in the method that may depend on $epsi$ is the \emph{number} of levels in the telescopic hierarchy. For systems in which the spectrum of the collision operator fall apart into a set of clearly separated clusters (each corresponding to a specific time scale), the number of levels equals the number of spectral clusters. In this situation, the computational cost is completely independent of $\epsi$.  When the collision operator represents a continuum of time scales, the number of projective integration levels increases logarithmically with $\epsi$. 

In this paper, we construct and evaluate telescopic projective integration methods for nonlinear Boltzmann BGK and Boltzmann kinetic equations. The methods are of arbitrary order in time, fully explicit, and general (they do not exploit any particular form of the collision operator).
The remainder of this paper is structured as follows. In Section~\ref{sec:models}, we start by presenting the Boltzmann and BGK equations that will be the subject of our simulations. We describe the different projective and telescopic projective integration methods in detail in Section \ref{sec:time_integrator}. (The spatial and velocity discretizations are standard. To make the manuscript self-contained, we present the corresponding numerical methods in Appendices~\ref{sub:weno} and~\ref{sub:FastSpectral}.) We discuss in Section \ref{sec:linearOpSpecProp} the spectral properties of the linearized collision operators, which will guide the choice of the method parameters, ensuring stability of the time integrators. Some numerical experiments are done in Section \ref{sec:results} to verify the theory developed in the two previous sections. We conclude in Section \ref{sec:conclusions}.

\section{Model equations} \label{sec:models}
In this article, we are interested in rarefied, collisional gases, and then we shall consider Boltzmann-like, collisional kinetic equations. We refer the reader to the classical works \cite{CIP,Villani:2002handbook} and the references therein for a more detailed introduction on this vast topic. 
For a given non-negative initial condition $f_0$, we will study a particle distribution function $\fepsi = \fepsi(\x,\v, t)$, for $t \geq 0$, $\x \in \Omega \subset \mathbb{R}^{D_x}$ and $\v \in \mathbb{R}^{D_v}$, solution to the initial-boundary value problem
\begin{equation} \label{eqCollision} 
	\left\{
	\begin{aligned}
		& \frac{\partial \fepsi}{\partial t} + \v \cdot \nabla_\x \fepsi = \frac{1}{\epsi} \collOp(\fepsi), \\
		&\, \\
		& \fepsi(\x, \v, 0) = f_{0}(\x,\v).
	\end{aligned}
	\right.
\end{equation}
The left hand side of equation \eqref{eqCollision} corresponds to a linear transport operator that comprises the convection of particles in space, whereas the right hand side contains the  collision operator that entails velocity changes due to particle collisions.  We postpone the description of boundary conditions until Section~\ref{sec:results}, where we discuss the numerical results.

We assume that the collision operator fulfils the following three assumptions:
\begin{enumerate}[label=\textbf{(H$\bm{_\arabic{*}}$)}, ref=\textbf{(H$\bm{_\arabic{*}}$)}]
	\item \label{hypConservations} Conservation of mass, momentum and kinetic energy:
	\begin{equation*}
		\int_{\R^{D_v}} \collOp(f)(\v) \, d\v = 0, \quad  \int_{\R^{D_v}} \collOp(f)(\v) \, \v \, d\v = \bm{0}_{\R^{D_v}}, \quad \int_{\R^{D_v}} \collOp(f)(\v) \, |\v|^2 \, d\v = 0;
	\end{equation*}
	 
	\item \label{hypEntropy} Dissipation of the Boltzmann entropy (H-theorem):
	\begin{equation*}
		\int_{\R^{D_v}} \collOp(f)(\v) \, \log(f)(\v) \, d\v \, \leq \, 0;
	\end{equation*}
			
	\item \label{hypEquilib} Its equilibria are given by Maxwellian distributions:
	\begin{equation*}
	    \collOp(f) \, = \, 0 \quad \Leftrightarrow \quad f = \cM_\v^{\rho, \vMacroMultiD, T} := \frac{\rho}{(2 \pi T)^{D_v/2}} \exp \left ( - \frac{|\v-\vMacroMultiD|^2}{2 T} \right ),
	\end{equation*}
    where the \emph{density} $\rho$, \emph{velocity} $\vMacroMultiD$ and	\emph{temperature} $T$  of the gas are computed from the distribution function $f$ as:
	\begin{equation} \label{eq:f_moments} 
		\rho = \int_{\R^{D_v}} f(\v)\,d\v, \quad \vMacroMultiD = \frac{1}{\rho}\int_{\R^{D_v}} \v f(\v) \, d\v, 
		\quad T = \frac{1}{D_v \rho} \int_{\R^{D_v}} \vert \vMacroMultiD - \v \vert^2 f(\v) \,d\v.
	\end{equation}
	      
\end{enumerate}  
Equation \eqref{eqCollision} with assumptions \ref{hypConservations}-\ref{hypEntropy}-\ref{hypEquilib}  describes numerous models such as the Boltzmann equation for elastic collisions \cite{Villani:2002handbook} or Fokker-Planck-Landau type equations\cite{alexandre2000entropy}. 
The parameter $\epsi > 0$ is the dimensionless Knudsen number, that is, the ratio between the mean free path of particles  and the length scale of observation. It determines the regime of the gas flow, for which we roughly identify the hydrodynamic regime $(\epsi \le 10^{-4})$, the transitional regime $(\epsi \in [10^{-4},10^{-1}])$, and the kinetic regime $(\epsi \ge 10^{-1})$.
Moreover, according to assumptions \ref{hypEntropy}-\ref{hypEquilib}, when $\epsi \to 0$, the distribution $\fepsi$ converges (at least formally) to a Maxwellian distribution, whose moments are solution to the compressible Euler system for perfect gases, given by:
\begin{equation} \label{eqHydroClosedEuler}
	\left\{
	\begin{aligned}
		& \partial_t \rho  + \diverg_\x  (\rho \, \vMacroMultiD ) = 0, 				  \\
		&\, \\
		& \partial_t(\rho \, \vMacroMultiD ) + \diverg_\x  \left(\rho \, \vMacroMultiD \otimes \vMacroMultiD \,+\, \rho \, T  \,{\rm\bf I}\right) \, =\, \bm{0}_{\R^{D_v}}, 				  \\
		&\, \\
		& \partial_t E + \diverg_\x \left ( \vMacroMultiD \left ( E +\rho \, T\right )  \right ) \,=\, 0,
	\end{aligned}
	\right.
\end{equation}
in which $E$ is the second moment of $\fepsi$, namely the total energy of the gas:
\begin{displaymath}
	E = \int_{\R^{D_v}} \vert \vMacroMultiD \vert^2 f(\v) \,d\v.
\end{displaymath}

In the following, we will present the two main collisional kinetic equations that we will consider in the remainder of this paper: the Boltzmann equation (Section~\ref{subsec:boltzmann_equation}) and the BGK equation (Section~\ref{subsec:bgk_equation}).
			
\subsection{Boltzmann equation} \label{subsec:boltzmann_equation}
 The Boltzmann equation constitutes the cornerstone of the kinetic theory of rarefied gases \cite{Villani:2002handbook, CIP}. In a dimensionless, scalar setting, it describes the evolution of the one-particle mass distribution function $\fepsi(\x,\v,t) \in \R^{+}$, solution to the model equation \eqref{eqCollision}, in which we still need to specify the collision operator $\collOper{\fepsi}$. 
The Boltzmann collision operator models binary elastic collisions between particles having pre-collisional velocities $(\v',\vstarb')$ and post-collisional velocities $(\v,\vstarb)$. In a two-dimensional velocity space, the pre- and post-collisional velocities are linked through the following parametrization:
\begin{equation*}
	\v' = \dfrac{\v + \vstarb}{2} + \dfrac{\abs{\v - \vstarb}}{2}\sigmavec, \qquad\quad
	\vstarb' = \dfrac{\v + \vstarb}{2} - \dfrac{\abs{\v - \vstarb}}{2}\sigmavec,
\end{equation*}
where $\sigmavec$ is the unit vector on the unit circle $\S^1 = \{\sigmavec \in \R^2: \abs{\sigmavec} = 1\}$ directed along the pre-collisional relative velocity $\vrel' = \v' - \vstarb'$:
\begin{equation*}
	\sigmavec = \frac{\vrel'}{\abs{\vrel'}} = \frac{\v' - \vstarb'}{\abs{\v' - \vstarb'}}.
\end{equation*}
The Boltzmann collision operator then reads:
\begin{align}
	\collOper{\fepsi} &= \int_{\R^2}\int_{\S^1} B(\abs{\vrel},\sigmavec)(f'f_*' - ff_*) d\sigmavec d\vstarb \notag \\
	&= \int_{\R^2}\int_{0}^{2\pi} B(\abs{\v-\v_*},\ths)(f'f_*' - ff_*) d\ths d\vstarb, \label{eq:collision_operator}
\end{align}
where  $\ths$ is the angle between $\vrel'$ and $\sigmavec$ and we used the shorthand notations $f = \fepsi(\v)$, $f_* = \fepsi(\v_*)$, $f^{'} = \fepsi(\v')$, and $f_*^{'} = \fepsi(\v_* ^{'})$. Furthermore, the non-negative function $B(\abs{\vrel},\sigmavec) \equiv B(\abs{\vrel},\ths)$ is the collision kernel, which, by physical arguments of invariance, only depends on the relative speed $\abs{\vrel} = \abs{\v-\v_*}$ and $\cos(\ths) = \vrel'/\abs{\vrel} \cdot \sigmavec$. The collision kernel $B$ contains all relevant microscopic information such as the kind of particles and type of interactions.
For instance, when particles interact via an inverse power law potential $\pot(r) = r^{-k+1}$ $(k > 2)$, with $r$ the inter-particle distance, $B$ factors as:
\begin{equation} \label{eq:collision_kernel_factors_sigma}
	B(\abs{\vrel},\ths) = \abs{\vrel}^\gamma b_\gamma(\ths), \qquad\quad
	\gamma = \frac{k-3}{k-1}.
\end{equation}
Notably, in the special case $k = 3$, the collision kernel in \eqref{eq:collision_kernel_factors_sigma} is independent of the relative speed $\abs{\vrel}$ and the resulting particles are known as Maxwellian particles. If, in addition, $b_0(\ths) = b_0$ is assumed constant, the particles are referred to as pseudo-Maxwellian particles.
In general (except for hard sphere and pseudo-Maxwellian particles), the angular collision kernel $b_\gamma$ in \eqref{eq:collision_kernel_factors_sigma} is expressed implicitly and contains a singularity for grazing collisions $(\ths \to 0)$, and its mathematical analysis, as well as its numerical simulations, can be very difficult \cite{alexandre2000entropy}. For that reason, the angular collision kernel is usually replaced by an integrable function by cutting off such grazing collision angles (Grad's cut-off assumption)  \cite{Cercignani1988}.

It is instructive to split the collision operator \eqref{eq:collision_operator} into a gain and loss operator as:
\begin{equation} \label{eq:collision_operator_gain_loss}
	\collOper{\fepsi} = \collOperPlus{\fepsi} - \collFreq(\fepsi)\fepsi(\v),
\end{equation}
where the gain operator is given by:
\begin{equation} \label{eq:collision_operator_gain}
	\collOperPlus{\fepsi} = \int_{\R^2}\int_{0}^{2\pi} B(\abs{\vrel},\ths)f'f_*' d\ths d\vstarb.
\end{equation}
The loss operator $\collOpMin(\fepsi) = \collFreq(\fepsi)\fepsi$ contains the collision frequency $\collFreq(\fepsi) \in \R^+$, which is defined by:
\begin{equation} \label{eq:collision_frequency_Boltzmann}
	\collFreq(\fepsi) = \int_{\R^2}\int_{0}^{2\pi} B(\abs{\vrel},\ths)f_* d\ths d\vstarb.
\end{equation}
(Evidently, equation \eqref{eq:collision_operator_gain_loss} is only valid if both integrals \eqref{eq:collision_operator_gain}-\eqref{eq:collision_frequency_Boltzmann} are convergent, which is certainly true for a cut-off collision kernel.)
In general, the collision frequency depends on the dimension of velocity space and the type of microscopic collisions. In particular, when considering pseudo-Maxwellian particles for ${D_v = 2}$, implying that $\gamma = 0$ and $b_0(\ths) = b_0$ constant, the collision kernel becomes much simpler and is given by ${B(\abs{\vrel},\ths) = b_0}$, that is, independent of $\abs{\vrel}$ and $\ths$.
In that case, the collision frequency \eqref{eq:collision_frequency_Boltzmann} can be explicitly computed as:
\begin{equation*}
	\collFreq(\fepsi) = b_0 \int_{0}^{2\pi} d\ths \int_{\R^2} f_* d\vstarb = 2\pi b_0\rho,
\end{equation*}
and the collision operator in \eqref{eq:collision_operator_gain_loss} reads:
\begin{equation} \label{eq:collision_operator_maxwellian_particles}
	\collOper{\fepsi} = \collOperPlus{\fepsi} - 2\pi b_0\rho \fepsi.
\end{equation}
We note that, for inverse power law potentials with $\gamma \neq 0$ (including hard sphere particles), the collision frequency in general depends on the density, temperature and collision kernel, see \cite{Struchtrup2005}.
Finally, it is easy (see \cite{CIP}) to check that the Boltzmann collision operator satisfies the hypotheses \ref{hypConservations}-\ref{hypEntropy}-\ref{hypEquilib}.

\subsection{BGK equation} \label{subsec:bgk_equation}
Due to the high-dimensional and complicated structure of the Boltzmann collision operator, the Boltzmann collision kernel~\eqref{eq:collision_operator} is often replaced by simpler collision models that capture most of its essential features. The most well-known approximation is the BGK model \cite{Bhatnagar1954}, which models collisions as a relaxation towards thermodynamic equilibrium. (Because of the moments-dependency of the equilibrium, this is still a nonlinear operator.) It is the most well regarded simplified model of the Boltzmann equation, and is almost universally used in the physics and numerics communities (see \cite{DimarcoPareschi2014} and the references therein for details). It is given by:
\begin{equation} \label{eq:bgk_equation}
	\partial_t \fepsi + \v \cdot \nabla_{\x} \fepsi = \frac{\nu}{\epsi}(\MaxBGKMultiD - \fepsi),
\end{equation}
where $\MaxBGKMultiD$ denotes the local Maxwellian distribution, which, for a $D_v$-dimensional velocity space, is given by:
\begin{equation} \label{eq:maxwellian} 
	\MaxBGKMultiD = \frac{\rho}{(2\pi T)^{D_v/2}} \exp{\left(-\frac{|\v-\vMacroMultiD|^2}{2T}\right)} := \cM_\v^{\rho,\vMacroMultiD,T}.
\end{equation}
Furthermore, in the BGK equation \eqref{eq:bgk_equation}, the collision frequency $\collFreq \in \R^+$ is in general derived from the Boltzmann collision operator and its expression is given in \eqref{eq:collision_frequency_Boltzmann}, while in the linearized setting, $\collFreq$ is independent of $\fepsi$ and is formulated in \eqref{eq:collision_frequency_linearized}, see Section~\ref{subsec:lin_boltzmann_equation}. Notably, when setting $\collFreq = \rho$ for $D_v = 2$, the BGK model matches the loss term of the Boltzmann collision operator for pseudo-Maxwellian particles close to equilibrium, see equation \eqref{eq:collision_operator_maxwellian_particles}.
The BGK collision operator satisfies by construction the hypotheses \ref{hypConservations}-\ref{hypEntropy}-\ref{hypEquilib}.

\section{Numerical method} \label{sec:time_integrator}

Now that we have introduced the model problems, we turn to the description of the numerical method that will be the focus of this paper.  Equation~\eqref{eqCollision} needs to be discretized in space, velocity and time. 

We discretize equation \eqref{eqCollision} in space using finite differences on a uniform, constant in time, periodic mesh with spacing $\dx$, consisting of $I$ mesh points $x_i=i\dx$, ${1 \le i \le I}$, with $I\dx=1$. In the numerical experiments of section~\ref{sec:results}, we use the classical WENO scheme \cite{shu:1998}, which we briefly recall in Appendix~\ref{sub:weno}. Next, we discretize velocity space by choosing $J$ discrete components denoted by $\v_j$. For the Boltzmann equation, that is, equation \eqref{eqCollision} with collision operator~\eqref{eq:collision_operator}, we use the fast spectral discretization of the Boltzmann operator, taken from \cite{MoPa:2006}. This method is recalled in Appendix~\ref{sub:FastSpectral}.

The semidiscrete numerical solution on this mesh is denoted by $f_{i,j}(t)$, where we have dropped the superscript $\epsi$ on discretized quantities. We then obtain a semidiscrete system of ODEs of the form:
\begin{equation}\label{eq:semidiscrete} 
	\dot{\fSDSys} = \Dtime(\fSDSys),  \qquad
	\Dtime(\fSDSys) = -\Dxv(\fSDSys) + \frac{1}{\epsi}\collOp\left (\fSDSys\right),
\end{equation}
where $\Dxv(\cdot)$ represents the finite difference discretization of the convective derivative $\v \cdot \nabla_\x$, and $\fSDSys$ is a vector of size $I \cdot J$.

In the remainder of this section, we describe the time discretization of the semi-discretized system~\eqref{eq:semidiscrete}, which is the novel element in the full discretization of equation~\eqref{eqCollision}.  We start in Section~\ref{sub:projective_integration} with the projective integration method, which aims at efficiently simulating systems with \emph{exactly two} time scales (one fast and one slow). In Section~\ref{sub:telescopic_integration}, we present the generalized telescopic projective integration method, which can deal with multiple fast time scales.

\subsection{Projective integration} \label{sub:projective_integration}

\begin{figure}[t]
	\begin{center}
		\inputtikz{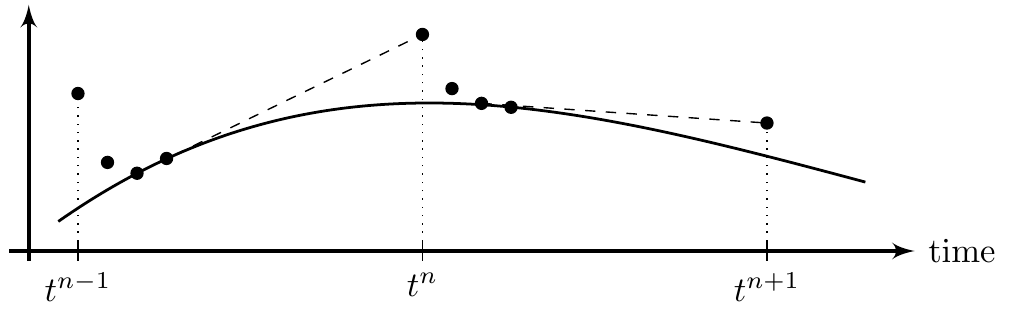}
	\end{center}
  	\vspace*{-0.5cm}\caption{\label{fig:proj_int} Sketch of projective integration. At each time, an explicit method is applied over a number of small time steps (black dots) so as to stably integrate the fast modes. As soon as these modes are sufficiently damped, the solution is extrapolated using a much larger time step (dashed lines). }
\end{figure}

Projective integration \cite{Gear2003projective} is a method that is tailored to problems with exactly two distinct time scales. As such, in the context of kinetic equations, it matches nicely with the spectral properties of a linear BGK equation, as was shown in~\cite{Lafitte2012}.  Projective integration combines a few small time steps with a naive (\emph{inner}) timestepping method (here, a direct forward Euler discretization) with a much larger (\emph{projective, outer}) time step. The idea is sketched in figure \ref{fig:proj_int}.

\paragraph{Inner integrators.} 
At the innermost leve, we introduce a uniform time mesh with time step $\dt$ and discrete time instants $t^k=k\dt$. At this leve, we choose the (explicit) forward Euler method with time step $\dt$, for which we will, later on, use the shorthand notation:
\begin{equation} \label{eq:fe_scheme} 
	\fSDSys^{k+1} = \Sdt(\fSDSys^{k}) = \fSDSys^k + \dt\Dtime(\fSDSys^k), \qquad k = 0, 1, \ldots.
\end{equation}
The purpose of the inner integrator is to capture the fastest components in the numerical solution of system \eqref{eq:semidiscrete} and to sufficiently damp these out. We only require the innermost integrator to be stable for these components. The size of the inner time step~$\dt$ and the required number of inner steps~$K$ will depend on the spectral properties of the semidiscretization~\eqref{eq:semidiscrete}. This will be studied in Section~\ref{sec:linearOpSpecProp}.

\paragraph{Outer integrators.} In system \eqref{eq:semidiscrete}, the small parameter $\epsi$ leads to a classical time step restriction of the form $\dt = O(\epsi)$ for the inner integrator. However, as $\epsi$ goes to $0$, we obtain the limiting system \eqref{eqHydroClosedEuler}, for which a standard finite volume/forward Euler method only needs to satisfy a CFL stability restriction of the form $\Dt \le C\dx$, with $C$ a constant that depends on the specific choice of the scheme.

In \cite{Lafitte2012}, it was proposed to use a projective integration method to accelerate such a brute-force integration; the idea, originating from \cite{Gear2003projective}, is the following. Starting from a computed numerical solution $\fSDSys^n$ at time $t^n=n\Dt$, one first takes $K+1$ \emph{inner} steps of size $\dt$ using~\eqref{eq:fe_scheme}, denoted as 
$\fSDSys^{n,k+1}$,
in which the superscripts $(n,k)$ denote the numerical solution at time ${t^{n,k}=n\Dt +k\dt}$. The aim is to obtain a discrete derivative to be used in the \emph{outer} step to compute $\fSDSys^{n+1} = \fSDSys^{n+1,0}$ via extrapolation in time:
\begin{align*}
	\fSDSys^{n+1} & = \fSDSys^{n,K+1} + (\Dt - (K + 1)\dt)\frac{\fSDSys^{n,K+1} - \fSDSys^{n,K}}{\dt}, \\
	 & = \fSDSys^{n,K+1} + M\dt\frac{\fSDSys^{n,K+1} - \fSDSys^{n,K}}{\dt},
\end{align*}
where $M = \Dt/\dt-(K+1)$.
Also the size of the (macroscopic) extrapolation step~$\Dt$ will result from the spectral analysis of the semidiscretization~\eqref{eq:semidiscrete} in section~\ref{sec:linearOpSpecProp}.

Higher-order projective Runge-Kutta (PRK) methods have been constructed \cite{LafitteLejonSamaey2015,LafitteMelisSamaey2017} by replacing each time derivative evaluation $\mathbf{k}_s$ in a classical Runge-Kutta method by $K+1$ steps of an inner integrator as follows:
\begin{align*}
	s = 1 :\;\; & 
	\begin{dcases} 
		\fSDSys^{n,k+1} &= \fSDSys^{n,k} + \dt\Dtime(\fSDSys^{n,k}), \qquad 0 \le k \le K \\ 
		\mathbf{k}_1 &= \dfrac{\fSDSys^{n,K+1} - \fSDSys^{n,K}}{\dt}
	\end{dcases} \\
	2 \le s \le S :\;\; & 
	\begin{dcases} 
		\fSDSys^{n+c_s,0}_s &= \fSDSys^{n,K+1} + (c_s\Dt-(K+1)\dt) \sum_{l=1}^{s-1}\dfrac{a_{s,l}}{c_s} \mathbf{k}_l, \\
		\fSDSys^{n+c_s,k+1}_s &= \fSDSys^{n+c_s,k}_s + \dt\Dtime(\fSDSys^{n+c_s,k}_s), \qquad 0 \le k \le K \\
		\mathbf{k}_s &= \dfrac{\fSDSys^{n+c_s,K+1}_s - \fSDSys^{n+c_s,K}_s}{\dt}
	\end{dcases} \\
	& \fSDSys^{n+1} = \fSDSys^{n,K+1} + (\Dt-(K+1)\dt)\sum_{s=1}^{S}b_s \mathbf{k}_s.
\end{align*}
To ensure consistency, the Runge-Kutta matrix $\mathbf{a}=(a_{s,i})_{s,i=1}^S$, weights ${\mathbf{b}=(b_s)_{s=1}^S}$, and nodes $\mathbf{c}=(c_s)_{s=1}^S$ satisfy the conditions $0\le b_s \le 1$ and $0 \le c_s \le 1,$ as well as:
\begin{equation*}
	\sum_{s=1}^Sb_s=1, \qquad \sum_{i=1}^{S-1} a_{s,i} =c_s, \quad 1 \le s \le S. 
\end{equation*}

\subsection{Telescopic projective integration} \label{sub:telescopic_integration}
In general, the stiff semidiscrete system \eqref{eq:semidiscrete}, contains more than two distinct time scales. In this section, we therefore describe an extension of projective integration, called \emph{telescopic projective integration} (TPI) and introduced in \cite{Gear2003telescopic}, that can handle multiple time scales.  This method has been studied in the context of linear BGK equations with multiple relaxation times in \cite{MelisSamaey2017}. 

\begin{figure}[t]
	\begin{center}
		\inputtikz{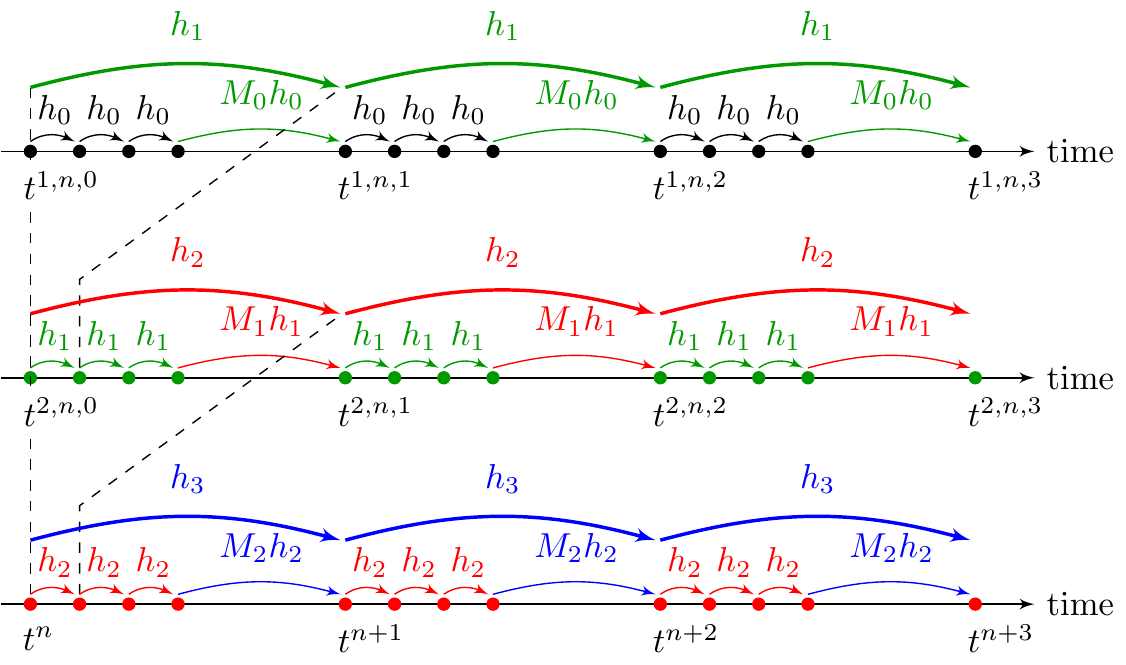}
	\end{center}
	\caption{\label{fig:tpi_sketch} A level-3 TPI method drawn for three outermost time steps $h_3$ (bottom row) with $K=2$ on all projective levels. The dots correspond to different time points at which the numerical solution is calculated. The time step and projective step size of each level $\ell=0,\ldots,2$ are denoted by $h_\ell$ and $M_\ell$, respectively. }
\end{figure}

Telescopic projective integration employs a number of projective integrator levels, which, starting from a base (\emph{innermost}) integrator, are wrapped around the previous level integrator \cite{Gear2003telescopic}. In this way, a hierarchy of projective integrators is formed in which each level (except the innermost and outermost one) fulfils both an inner and outer integrator role. This generalizes the idea of projective integration, which contains only one projective level wrapped around an inner integrator. 
The idea of a level-3 TPI method with $K=2$ on each projective level is sketched in figure \ref{fig:tpi_sketch}.
The different level integrators in a TPI method can in principle be selected independently from each other, but in general one selects a first order explicit scheme (the forward Euler scheme) for all but the outermost integrator level, whose order is chosen to meet the accuracy requirements dictated by the problem.
%

\paragraph{Innermost integrator}
We intend to integrate the semidiscrete system of equations \eqref{eq:semidiscrete} using a uniform time mesh with time step $h_0$ and discrete time instants $t^k=kh_0$. The innermost integrator of the TPI method is chosen to be the forward Euler (FE) method,
\begin{equation*}
	\fSDSys^{k+1} = \fSDSys^k + h_0\Dtime(\fSDSys^k). 
\end{equation*}
In the sequel, we use the following shorthand notation:
\begin{equation*} 
	\fSDSys^{k+1} = S_{0}(\fSDSys^{k}) \qquad (k = 0, 1, \ldots),
\end{equation*}
in which $S_0$ denotes the time stepper with corresponding time step $h_0$.
Also in the telescopic projective integration method, the purpose of the innermost integrator is only to capture the fastest components in the numerical solution of system~\eqref{eq:semidiscrete} and to sufficiently damp these out. As a consequence, it is ill-advised to use higher-order methods for the innermost integrator, see~\cite{MelisSamaey2017} for a more detailed discussion.

\paragraph{Projective (outer) levels} 
The telescopic projective integration method employs in general $L$ nested levels of  projective integration that are constructed around the innermost integrator. In \cite{Gear2003telescopic}, the method has been introduced in a recursive way.  Here, following~\cite{MelisSamaey2017}, we describe the method in an alternative way, to make the presentation more similar to that of classical projective integration. 

To keep track of the time instant at which the numerical solution is computed throughout the TPI method and at the same time desiring a compact notation, in what follows, we employ superscript triplets of the form $(\ell,n,k_\ell)$ where $\ell$ denotes the integrator level ranging from $0$ (innermost) to $L-1$, $n$ represents the index of the current outermost integrator time $t^n=nh_L$, and $k_\ell$ corresponds to the iteration index of the integrator on level $\ell$. The numerical time on each level $\ell = 0, \ldots, L-1$ is then defined as (see also figure \ref{fig:tpi_sketch}):
\begin{equation} \label{eq:tpi_time}
	t^{\ell,n,k_\ell} = nh_L + \sum_{\ell'=\ell}^{L-1} k_{\ell'}h_{\ell'}.
\end{equation}
Notice that, for a certain level $\ell$, this time requires the iteration indices $k_{\ell'}$ of all its outer integrators. Therefore, it incorporates a memory that keeps up with the current time instants at which the outer integrators of a given level $\ell$ integrator have arrived at and is necessary to take into account to correctly reflect the numerical time of the solution on each level $\ell$.

Starting from a computed numerical solution $\fSDSys^n$ at time $t^n=nh_L$, one first takes $K_0+1$ steps of size $h_0$ with the innermost integrator:
\begin{equation} \label{eq:semidiscrete_innermost_scheme}
	\fSDSys^{0,n,k_0+1} =S_{0}(\fSDSys^{0,n,k_0}) \qquad (0 \le k_0 \le K_0),
\end{equation}
in which $\fSDSys^{0,n,k_0}$ corresponds to the numerical solution at time $t^{0,n,k_0}$ calculated by the innermost integrator. Since all outer integrator iteration indices $k_{\ell'}$, $\ell' = 1,\ldots,L-1$ are initially zero in \eqref{eq:tpi_time}, we have ${t^{0,n,k_0} = nh_L+k_0h_0}$. The repeated action \eqref{eq:semidiscrete_innermost_scheme} of the innermost integrator is depicted by small black arrows in the upper row of figure \ref{fig:tpi_sketch}, for which we chose $K_0=2$.

In the telescopic projective integration framework, the scheme is set up from the lowest level up to the highest level. The aim is to obtain a discrete derivative to be used on each level to eventually compute $\fSDSys^{n+1} = \fSDSys^{0,n+1,0}$ via extrapolation in time. Using the innermost integrator iterations \eqref{eq:semidiscrete_innermost_scheme}, we perform the extrapolation by a projective integrator on level 1, written as:
\begin{equation} \label{eq:tpfe_first_level} 
	\fSDSys^{1,n,1} = \fSDSys^{0,n,K_{0}+1} + \left(M_{0}h_{0}\right)\frac{\fSDSys^{0,n,K_{0}+1} - \fSDSys^{0,n,K_{0}}}{h_{0}},
\end{equation}
which corresponds to the projective forward Euler (PFE) method. In \eqref{eq:tpfe_first_level}, $\fSDSys^{1,n,1}$ represents the numerical solution at time $t^{1,n,1}$ calculated by one iteration of the first level projective integrator. Since $k_1=1$ and all its outer integrator iteration indices $k_{\ell'}$, ${\ell'=2,\ldots,L-1}$ are still zero in \eqref{eq:tpi_time}, we have ${t^{1,n,1} = nh_L + h_1}$. One step of the first level integrator is visualized by a large green arrow in the upper row of figure \ref{fig:tpi_sketch}.
By repeating this idea, we construct a hierarchy of projective integrators on levels $\ell=1,\ldots,L-1$, given by:
\begin{equation} \label{eq:tpfe_projective_levels} 
	\fSDSys^{\ell,n,k_\ell+1} = \fSDSys^{\ell-1,n,K_{\ell-1}+1} + \left(M_{\ell-1}h_{\ell-1}\right)\frac{\fSDSys^{\ell-1,n,K_{\ell-1}+1} - \fSDSys^{\ell-1,n,K_{\ell-1}}}{h_{\ell-1}},
\end{equation}
in which, on each level $\ell$, we iterate over $k_{\ell} = 0, \ldots, K_{\ell}$.
In \eqref{eq:tpfe_projective_levels}, $\fSDSys^{\ell,n,k_\ell}$ denotes the numerical solution at time $t^{\ell,n,k_\ell}$ calculated by the projective integrator on level $\ell$. According to \eqref{eq:tpi_time}, this time depends on the values $k_{\ell'}$, $\ell'=\ell+1,\ldots,L-1$ of all of its outer integrators.
In figure \ref{fig:tpi_sketch}, these projective integrator steps are shown by long arrows for each level $\ell=1,\ldots,3$. Ultimately, the outermost integrator on level $L$ computes $\fSDSys^{n+1}$ as:
\begin{equation} \label{eq:tpfe_outermost}
	\fSDSys^{n+1} = \fSDSys^{L-1,n,K_{L-1}+1} + (M_{L-1}h_{L-1})\frac{\fSDSys^{L-1,n,K_{L-1}+1} - \fSDSys^{L-1,n,K_{L-1}}}{h_{L-1}}.
\end{equation}
Since the outermost integrator \eqref{eq:tpfe_outermost} also constitutes a PFE scheme, the telescopic method resulting from the hierarchy of projective levels \eqref{eq:tpfe_projective_levels}-\eqref{eq:tpfe_outermost} is called telescopic projective forward Euler (TPFE).

It is straightforward to implement higher-order extensions of the outermost integrator, as is done in~\cite{MelisSamaey2017}.  We mention the projective Runge-Kutta methods of order 2 and 4, leading to TPRK2 and TPRK4 method  in the telescopic case. In general, the outermost integrator in a TPRK method replaces each time derivative evaluation $\mathbf{k}_s$ in a classical Runge-Kutta method by $K_{L-1}+1$ steps of its inner integrator on level $L-1$. Using \eqref{eq:tpfe_projective_levels} with $\ell=L-1$, the first stage in a TPRK method calculates the time derivative $\mathbf{k}_1$ as:
\begin{equation} \label{eq:tprk_stage_1}
	\mathbf{k}_1 = \dfrac{\fSDSys^{L-1,n,K_{L-1}+1} - \fSDSys^{L-1,n,K_{L-1}}}{h_{L-1}}.
\end{equation}
Computing $\mathbf{k}_s$ on any other stage $s \ge 2$ requires evaluating time derivatives at the intermediate times ${t^{n+c_s} = (n+c_s)h_L}$. Similarly to \eqref{eq:tprk_stage_1}, these are computed as:
\begin{equation} \label{eq:tprk_stage_s}
	\mathbf{k}_s = \dfrac{\fSDSys^{L-1,n+c_s,K_{L-1}+1} - \fSDSys^{L-1,n+c_s,K_{L-1}}}{h_{L-1}}.
\end{equation}
However, since the numerical solution at time $t^{n+c_s}$ in equation \eqref{eq:tprk_stage_s} is not available, we use the integrator on level $L-1$ to approximate it as follows:
\begin{equation*}
	\begin{dcases}
		\fSDSys^{L-1,n+c_s,0} = \fSDSys^{L-1,n,K_{L-1}+1} + (c_sh_L-(K_{L-1}+1)h_{L-1}) \sum_{i=1}^{s-1}\dfrac{a_{s,i}}{c_s} \mathbf{k}_i \\
		\fSDSys^{L-1,n+c_s,k_{L-1}+1} = \fSDSys^{L-2,n+c_s,K_{L-2}+1} \\
		\hspace{3.1cm} + \,\left(M_{L-2}h_{L-2}\right)\frac{\fSDSys^{L-2,n+c_s,K_{L-2}+1} - \fSDSys^{L-2,n+c_s,K_{L-2}}}{h_{L-2}},
	\end{dcases}
\end{equation*}
in which the second equation iterates over $0 \le k_{L-1} \le K_{L-1}$.
Ultimately, the outermost integrator of a TPRK method is written as:
\begin{equation*}
	\fSDSys^{n+1} = \fSDSys^{L-1,n,K_{L-1}+1} + (M_{L-1}h_{L-1})\sum_{s=1}^{S}b_s \mathbf{k}_s.
\end{equation*}
The consistency conditions on the Runge-Kutta matrix $\mathbf{a}=(a_{s,i})_{s,i=1}^S$, weights $\mathbf{b}=(b_s)_{s=1}^S$, and nodes ${\mathbf{c}=(c_s)_{s=1}^S}$ are still valid in this setting \cite{MelisSamaey2017}.
In the numerical experiments in Section~\ref{sec:results}, we will use the projective Runge-Kutta method of order 4 as outermost integrator.

\section{On linearized operators and spectral properties}
\label{sec:linearOpSpecProp}
Telescopic projective integration methods are very versatile, but require choosing a relatively large number of method parameters: the size of the time steps $h_\ell$ at each level, the number $K_\ell$ of inner steps at each level, as well as the number $L$ of telescopic levels. As these choices are dictated mainly by stability requirements, they crucially depend on the spectrum of the collision operator. The analysis of this spectrum for the problems of Section~\ref{sec:models} is the focus of this Section.
 In \cite{Lafitte2012, LafitteLejonSamaey2015}, the spectrum of the collision operator was analysed for linear kinetic equations in the diffusive and hyperbolic scalings. An extension to kinetic relaxations of a nonlinear hyperbolic conservation law was presented in \cite{LafitteMelisSamaey2017}.  In all of these settings, the spectrum of the collision operator turned out to consist of exactly two well-separated time scales, and projective integration was therefore sufficient. A first study of telescopic projective integration for linear kinetic equations with multiple relaxation times was presented in \cite{MelisSamaey2017}.

In this section, we devise a general framework of linear operators in which linearizations of both the BGK and Boltzmann equation can be studied. This framework will allow determining suitable method parameters for the (telescopic) projective integration of the Boltzmann and nonlinear BGK equations. In addition, it allows embedding the linear kinetic equations that were studied in \cite{LafitteMelisSamaey2017,MelisSamaey2017}. For the reader's convenience, we restrict the exposition to the case when $D_v = 2$. However, all the results of this section can be extended straigthforwardly to the $D_v=3$ case, at the cost of heavier notations. In Section~\ref{subsec:lin_bgk}, we build this framework for BGK equations and draw conclusions on their spectral properties.  Afterwards, in Section~\ref{subsec:lin_boltzmann_equation}, we extend this framework to include the Boltzmann equation. We discuss the selection of suitable method parameters for (telescopic) projective integration in Section~\ref{sec:method_param}.

\subsection{Linearized BGK models and their spectra} \label{subsec:lin_bgk}
We first recall the linearization of the BGK equation~\eqref{eq:bgk_equation}, as it has been described in \cite{Cercignani1988}. We then show in section~\ref{sec:linear_kinetic} how the simpler linear kinetic equations that were analyzed in \cite{LafitteMelisSamaey2017,MelisSamaey2017} fit in this framework. Finally, in section~\ref{sec:bgk_spectrum}, we discuss how the analysis of the spectrum of linear kinetic equations in \cite{LafitteMelisSamaey2017,MelisSamaey2017} generalizes to the linearization of the full BGK equation. 

\subsubsection{Linearized BGK equation}\label{sec:linearized_bgk} 
In \cite{Cercignani1988}, it is shown that the linearized BGK operator can be formulated as:
\begin{equation} \label{eq:linearized_bgk_maxwellian}
	\cM_1(\fepsi)(\x,\v,t) = \sum_{k=0}^{D_v+1} \Psi_k(\v) \langle\Psi_k, \fepsi\rangle(\x,t),
\end{equation}
in which the scalar product is defined by:
\begin{equation} \label{eq:scalar_product_Hilbert}
	\langle g,h \rangle= \int_{\R^{D_v}} g(\v)\overline{h(\v)} \frac{1}{2\pi}\exp\left(\frac{-\abs{\v}^2}{2}\right) d\v.
\end{equation}
Furthermore, the basis functions $\Psi_k(\v)$ in \eqref{eq:linearized_bgk_maxwellian} represent an orthogonal basis for the space spanned by
\begin{equation*}
V:= \left\{1, v^{x_1}, \ldots, v^{x_{D_v}}, \frac{\abs{\v}^2}{2}\right\},
\end{equation*}
in which the superscript $x_d$ on $v$ indicates that $v^{x_d}$ is the component of $\v$ in the $d$-th velocity dimension.
The space $V$ corresponds to the set of elementary collision invariants $(\psi_k(\v))_{k=0}^{{D_v}+1}$. We construct an orthonormal basis, i.e., we seek basis functions $Psi_k$ such that:
\begin{equation} \label{eq:normalization_Psi}
	(\Psi_k, \Psi_j) = \delta_{kj}, \qquad k,j \in \{0,\ldots,{D_v}+1\},
\end{equation}
in which $\delta_{kj}$ denotes the Kronecker delta.
A straightforward application of the Gram-Schmidt process then allows to compute the desired set of orthonormal functions $\Psi_k(\v)$ satisfying \eqref{eq:normalization_Psi} as:
\begin{equation} \label{eq:Psi_normalized}
	\big(\Psi_0(\v), \ldots, \Psi_{D_v+1}(\v)\big) = \left(1, v^{x_1},\ldots, v^{x_{D_v}}, \frac{\abs{\v}^2 - 2}{2}\right).
\end{equation}
Ultimately, using the linearized Maxwellian \eqref{eq:linearized_bgk_maxwellian}, the linearized version of the full BGK equation \eqref{eq:bgk_equation} reads:
\begin{equation} \label{eq:linearized_bgk_equation_framework}
	\partial_t \fepsi + \v \cdot \nabla_{\x} \fepsi = -\frac{\collFreq}{\epsi}(\cI - \Pi_\text{BGK})\fepsi,
\end{equation}
where $\Pi_\text{BGK}$ is the following rank-${D_v}+2$ projection operator:
\begin{equation} \label{eq:projection_operator_bgk}
	\Pi_\text{BGK} \fepsi = \sum_{k=0}^{{D_v}+1} \Psi_k(\v)\langle\Psi_k, \fepsi\rangle.
\end{equation} 
 
\subsubsection{Linear kinetic models\label{sec:linear_kinetic}} 
 In \cite{LafitteMelisSamaey2017}, the spectrum of a specific class of linear, hyperbolically scaled, kinetic equations was studied. In a scalar, two-dimensional setting (where we shall set $x := x_1$ and $y := x_2$), the following equation was proposed:
\begin{equation} \label{eq:artificial_bgk_equation}
	\partial_t \fepsi + \v \cdot \nabla_{\x} \fepsi = \frac{1}{\epsi}(\MaxArtScal{\v} - \fepsi),
\end{equation}
with an artificial Maxwellian distribution given by:
\begin{equation} \label{eq:artificial_maxwellian}
	\MaxArtScal{\v} = \rho^\epsi (1+  v^x + v^y),
\end{equation}
in which $\rho^\epsi$ is linked to $\fepsi$ by averaging over the velocity space: 
\begin{equation*}
	\rho^\epsi = \int_{\R^2} \fepsi(\v) \frac{1}{2\pi} \exp\left(\frac{-\abs{\v}^2}{2}\right )d\v
\end{equation*}
see also equation~\eqref{eq:f_moments}.

In \cite{MelisSamaey2017}, multiple relaxation time linearized BGK equations of the following form are considered:
\begin{equation} \label{eq:artificial_linearized_bgk_equation}
	\partial_t \fepsi + \v \cdot \nabla_{\x} \fepsi = \frac{\collFreq}{\epsi}(\MaxBGKLinMultiD - \fepsi),
\end{equation}
with linearized Maxwellian given by:
\begin{equation} \label{eq:artificial_maxwellian_linearized} 
	\MaxBGKLinMultiD = \rho^\epsi(1 + v^x)(1 + v^y).
\end{equation}

To fit both of these model problems in the general framework of section~\ref{sec:linearized_bgk}, we rewrite the artificial Maxwellians in \eqref{eq:artificial_maxwellian} and \eqref{eq:artificial_maxwellian_linearized} as:
\begin{align}
	\cM_{-2}(\fepsi)(\x,\v,t) &= \Psi_{-2}(\v) \rho^\epsi(\x,t),& \quad \Psi_{-2}(\v) &= (1 + v^x + v^y), \label{eq:framework_artificial} \\
	\cM_{-1}(\fepsi)(\x,\v,t) &= \Psi_{-1}(\v) \rho^\epsi(\x,t),& \quad \Psi_{-1}(\v) &= (1 + v^x)(1 + v^y). \label{eq:framework_linearized_artificial}
\end{align}
With these choices of the Maxwellian distribution, we can summarize the linear kinetic models \eqref{eq:artificial_bgk_equation} and \eqref{eq:artificial_linearized_bgk_equation} as:
\begin{equation*}
	\partial_t \fepsi + \v \cdot \nabla_{\x} \fepsi = -\frac{\collFreq}{\epsi}(\cI - \Pi_k)\fepsi,
\end{equation*}
where $\cI$ is the identity operator and $\Pi_k$ is the following rank-1 projection operator:
\begin{equation} \label{eq:projection_operator_artificial}
	\Pi_k \fepsi = \Psi_k(\v) \int_{\R^2} \fepsi(\v) \frac{1}{2\pi}\exp\left(\frac{-\abs{\v}^2}{2}\right) d\v,
\end{equation}
with either $k = -2$ \eqref{eq:framework_artificial} or $k = -1$ \eqref{eq:framework_linearized_artificial}.

\subsubsection{Linearized BGK spectrum\label{sec:bgk_spectrum}} 
The above considerations indicate that the structure of the linearized Maxwellian \eqref{eq:linearized_bgk_maxwellian} and the linearized BGK projection operator \eqref{eq:projection_operator_bgk} are almost identical to those in \eqref{eq:framework_artificial}-\eqref{eq:framework_linearized_artificial} and \eqref{eq:projection_operator_artificial}, respectively. Indeed, the linear kinetic projection operators $\Pi_{-2}$ and $\Pi_{-1}$ in \eqref{eq:projection_operator_artificial} match the first three terms of $\Pi_\text{BGK}$ in \eqref{eq:projection_operator_bgk} and only differ in the last term. This can be seen by using the orthonormal set of basis functions \eqref{eq:Psi_normalized} and the scalar product \eqref{eq:scalar_product_Hilbert}, and subsequently rewriting $\Pi_{-2}$ and $\Pi_{-1}$ as:
\begin{equation} \label{eq:projection_operator_artificial_rewritten}
	\begin{aligned}
		\Pi_{-2} \fepsi &= \sum_{k=0}^2 \Psi_k(\v)(\Psi_k, \fepsi) \\
		\Pi_{-1} \fepsi &= \sum_{k=0}^2 \Psi_k(\v)(\Psi_k, \fepsi) + \Psi_1(\v)\Psi_2(\v)(\Psi_0, \fepsi).
	\end{aligned}
\end{equation}
We can thus view the linear kinetic models \eqref{eq:artificial_bgk_equation} and \eqref{eq:artificial_linearized_bgk_equation} used in \cite{LafitteMelisSamaey2017,MelisSamaey2017}, respectively, as a special simplified case of the linearized BGK equation.

Since the linearized BGK operator \eqref{eq:projection_operator_bgk} was shown to be nearly identical to the relaxation operators of the linear kinetic models in \eqref{eq:projection_operator_artificial_rewritten}, we expect the spectral properties of the linearized BGK equation \eqref{eq:linearized_bgk_equation_framework} to closely resemble those in \cite{LafitteMelisSamaey2017} (for $\collFreq = 1$) or \cite{MelisSamaey2017} (for $\collFreq = \rho$).
Therefore, it is expected that the construction of stable PI methods for the full BGK equation \eqref{eq:bgk_equation} with $\collFreq = 1$ and stable TPI methods for \eqref{eq:bgk_equation} with $\collFreq = \rho$ is practically identical to that in \cite{LafitteMelisSamaey2017} and \cite{MelisSamaey2017}, respectively.  The choice of method parameters for the full BGK equation will be discussed more closely in section~\ref{sec:method_param}.

\subsection{Linearized Boltzmann equation and its spectrum} \label{subsec:lin_boltzmann_equation}

\subsubsection{Linearization of the Boltzmann equation} \label{subsubsec:lin_boltzmann_equation}
To simplify the analysis of the Boltzmann equation \eqref{eqCollision}, it is customary to linearize the collision operator $\collOper{\fepsi}$ around the global Maxwellian distribution $\MaxGlobMultiD = \cM_\v^{1,0,1}$, which, for $D_v = 2$, is given by:
\begin{equation} \label{eq:global_maxwellian}
	\cM_\v^{1,0,1} = \frac{1}{2\pi} \exp\left(-\frac{\abs{\v}^2}{2}\right),
\end{equation}
see, for instance, \cite{Cercignani1988,Golse2005,Saint-Raymond2009}. Subsequently, we consider small fluctuations of $\fepsi$ around the global equilibrium \eqref{eq:global_maxwellian}, that is:
\begin{equation} \label{eq:linearization_f}
	\fepsi(\v) = \cM_\v^{1,0,1}\big(1 + \gepsi(\v)\big).
\end{equation}
Using a similar shorthand notation as before, that is, $g' = \gepsi(\v')$, $\cM' = \cM_{\v'}^{1,0,1}$, and so on, we linearize the quadratic terms in \eqref{eq:collision_operator} as:
\begin{equation} \label{eq:linearizing_ff_ff}
	f'f_*' - ff_* = \cM\cM_*\big(g' + g_*' - g - g_*\big),
\end{equation}
where we neglected second-order fluctuations in the second equality, and we used that $\cM'\cM_*' = \cM\cM_*$ for any Maxwellian distribution $\cM$. By substituting \eqref{eq:linearization_f} and \eqref{eq:linearizing_ff_ff} into the Boltzmann equation \eqref{eqCollision}, and exploiting that $\cM = \cM_\v^{1,0,1}$ depends only on velocity $\v$, we obtain the \emph{linearized Boltzmann equation} as:
\begin{equation} \label{eq:linearized_Boltzmann_equation}
	\partial_t \gepsi + \v \cdot \nabla_{\x} \gepsi = \frac{1}{\epsi}\collOpLin{\gepsi},
\end{equation}
where $\collOpLin{\gepsi}$ is the \emph{linearized Boltzmann collision operator}, which reads:
\begin{equation} \label{eq:linearized_collision_operator}
	\collOpLin{\gepsi} = \int_{\R^2}\int_{0}^{2\pi} B(\abs{\vrel},\ths)\cM_*(g_*' + g' - g_* - g) d\ths d\vstarb,
\end{equation}
see also \cite{CIP,Saint-Raymond2009}. Moreover, the operator $\collOpLin{\gepsi}$ can be cast in the following form \cite[IV.5]{Cercignani1988}:
\begin{equation} \label{eq:linearized_collision_operator_K_nu}
	\collOpLin{\gepsi}(\v) = \cK\gepsi(\v) - \collFreq(\abs{\v})\gepsi(\v),
\end{equation}
where $\collFreq(\abs{\v})$ is a local multiplication operator termed the \emph{linearized collision frequency} that depends only on the magnitude of $\v$ and is defined as:
\begin{equation} \label{eq:collision_frequency_linearized}
	\collFreq(\abs{\v}) = \int_{\R^2}\int_{0}^{2\pi} B(\abs{\vrel},\ths)\cM_* d\ths d\vstarb,
\end{equation}
and $\cK$ is a non-local integral operator containing the remaining three terms in \eqref{eq:linearized_collision_operator}. 

Using \eqref{eq:linearized_collision_operator_K_nu} and writing $\fepsi$ instead of $\gepsi$, we rewrite \eqref{eq:linearized_Boltzmann_equation} as:
\begin{equation} \label{eq:linearized_Boltzmann_equation_framework}
	\partial_t \fepsi + \v \cdot \nabla_\x \fepsi = -\frac{1}{\epsi}\big(\collFreq(\abs{\v})\cI - \cK\big)\fepsi.
\end{equation}
For inverse power law potentials under Grad's cut-off assumption and for hard sphere particles, it can be proven that $\cK$ is a compact operator on $L^2(\R^2)$, see \cite{Cercignani1988,Golse2005}. This implies that it maps the unit ball of $\R^2$ onto a finite-dimensional space \cite{Cercignani1988}. In that sense, it shares properties with the finite rank operators described previously.
However, the linearized operator $\collOpLin = \cK - \collFreq(\abs{\v})\cI$ does not have finite rank. Furthermore, one could write formally:
\begin{equation*}
	\cK = \collFreq(\abs{\v})(\Pi_\text{BGK} + \cR),
\end{equation*}
for a certain remainder operator $\cR$. Then, one can write \eqref{eq:linearized_Boltzmann_equation_framework} as:
\begin{equation*}
	\partial_t \fepsi + \v \cdot \nabla_\x \fepsi = -\frac{\collFreq}{\epsi} (\cI - \Pi_\text{BGK})\fepsi - \frac{\collFreq}{\epsi} \cR \fepsi,
\end{equation*}
which provides a connection with the linearized BGK projection operator.

%

\subsubsection{Linearized Boltzmann spectrum.}
Analyzing the spectrum of the linearized Boltzmann collision operator $\collOpLin = \cK - \collFreq(\abs{\v})\cI$ in \eqref{eq:linearized_Boltzmann_equation_framework} is more involved than in the BGK case. In general, the linearized Boltzmann collision operator has a spectrum that consists of (i) a non-empty essential (purely continuous) part that is entirely determined by the continuous spectrum of $-\collFreq(\abs{\v})\cI$, and (ii) a set of discrete eigenvalues that is influenced by the operator $\cK$, see, for instance, \cite{Cercignani1988} or \cite{BarangerMouhot2005,Ellis1975}. In contrast, the spectrum of the linear kinetic relaxation operators and linearized BGK operator only consists of discrete eigenvalues. However, for Maxwellian particles with angular cut-off and, in particular, for pseudo-Maxwellian particles, it is known that the spectrum of $\collOpLin$ contains only discrete eigenvalues spread inside the interval $[-\collFreq(0),0]$ with $\collFreq(\abs{\v})$ given in \eqref{eq:collision_frequency_linearized} \cite{Cercignani1988}. 

\begin{figure}
 \begin{center}
   \includegraphics[scale=1]{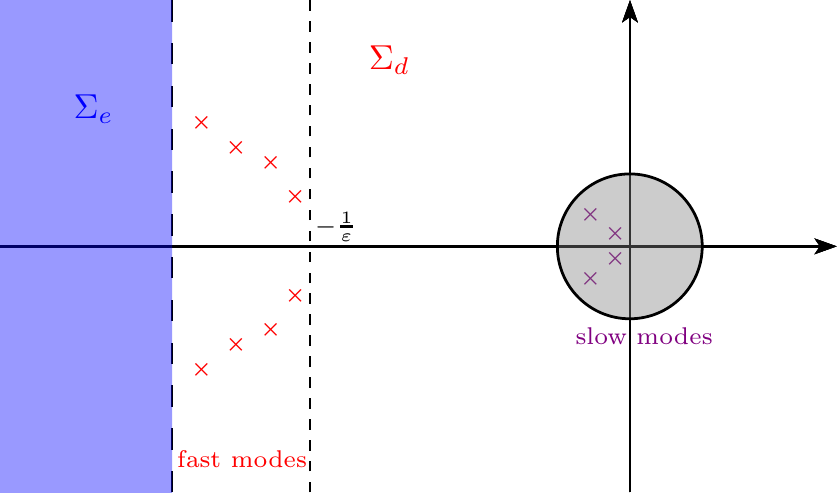}
 \end{center}
 \caption{\label{fig:spectrumBoltzFourier} Spectrum of the Fourier transformed linearized Boltzmann operator, for small radial frequencies.}
\end{figure}

We have seen in Section \ref{subsubsec:lin_boltzmann_equation} that the linearized version of the full Boltzmann equation reads:
\begin{equation*} 
	\partial_t \fepsi + \v \cdot \nabla_{\x} \fepsi = -\frac{1}{\epsi}(\nu(|\v|)\cI - \cK)\fepsi.
\end{equation*}
Let us apply the Fourier transform in the physical space: since the collision operator depends only on the velocity magnitude $\abs{\v}$, the only difference in the equation will be that the free transport term $\v \cdot \nabla_\x$ will become a multiplication operator (by $i \bm \gamma \cdot \v$, where $\bm \gamma$ is the spatial Fourier variable). 
One can then write the Fourier-transformed linear Boltzmann equation as:
\begin{equation} \label{eq:linearBoltzFourier}
	\partial_t \hepsi = \frac 1\epsi \cK \hepsi - \left ( \nu(|\v|)/\epsi + i \, \epsi \bm \gamma \cdot \v \right ) \hepsi,
\end{equation}
where $\hepsi$ is the Fourier transform in space of $\fepsi$. 
Hence, the evolution of $\hepsi$ is given by a compact perturbation of a (complex-valued) multiplication operator. 
It was proven in a series of papers that the spectrum of this Fourier-transformed collision operator has the following behavior as a function of $|\bm \gamma|$ and $\epsi$:

\begin{thm}[\cite{Nicolaenko71}, Section 2, and \cite{Ellis1975}, Theorem 3.1] \label{thm:linBoltzSpectrum}
	The spectrum of the right hand side of equation \eqref{eq:linearBoltzFourier} consists of an essential part $\Sigma_e$ located to the left of a vertical line of negative real part and a discrete spectrum $\Sigma_d$ composed of:
    \begin{itemize}
	    \item \emph{fast modes:} eigenvalues located at a distance at least $1/\epsi$ to the left of the imaginary axis;         
        \item \emph{slow modes:} if $|\epsi| \ll 1$, there are \emph{exactly} $4$ eigenvalues branches given by:
	    \[
	    \lambda^{(j)}(|\bm \gamma|) := i \, \lambda^{(j)}_1 \, \epsi |\bm \gamma| - \lambda^{(j)}_2 \, \epsi^2 |\bm \gamma|^2 + \mathcal O \left ( \epsi^3 |\bm \gamma|^3 \right ), \quad j \in \{0,1,2,3\},
	    \]
		for explicit constants $\lambda_1^{(j)} \in \R$ and $\lambda_2^{(j)}>0$.
	\end{itemize} 
\end{thm}
A sketch of this result can be found in figure \ref{fig:spectrumBoltzFourier}.

We observe in Theorem \ref{thm:linBoltzSpectrum} that the discrete eigenvalues form a fast and a slow cluster, justifying the use of projective integration with easily computable parameters, in the spirit of \cite{LafitteMelisSamaey2017,MelisSamaey2017}. Nevertheless, the presence of an essential spectrum is one of the reasons that one must use \emph{telescopic} projective integration to solve the full Boltzmann equation: this spectral decomposition will give rise to new clusters of eigenvalues at the discrete level. 

\begin{rem}
	Note that one can mimic the proof of this result for the simpler BGK operator with $\nu = 1$ to obtain the same spectral behavior of the linearized operator without the essential part, justifying at the continuous level the results from \cite{LafitteMelisSamaey2017}.
\end{rem}

\subsection{Method parameters for projective and telescopic projective integration} \label{sec:method_param}
It still remains to select appropriate parameter values for the projective or telescopic projective integration methods. These are determined by ensuring that all eigenvalues of the kinetic problem under study fall within the stability region of the full projective method.
In sections \ref{subsec:lin_bgk} and \ref{subsec:lin_boltzmann_equation}, we revealed that the spectra of the linearized kinetic equations either appear in two stationary eigenvalue disks (linearized BGK equation with $\collFreq = 1$) or are continuously spread along (a part of) the negative real axis (linearized BGK equation with $\collFreq = \rho$ and linearized Boltzmann equation).
Since the construction of stable projective methods for both well-separated as well as continuously spread spectra is studied in previous works \cite{LafitteMelisSamaey2017} and \cite{MelisSamaey2017}, respectively, we take over the main results here, which are summarized below.

\subsubsection{Stationary, well-separated spectrum}
For the linearized BGK equation with $\nu = 1$, which falls into the class of kinetic models studied in \cite{LafitteMelisSamaey2017}, it was shown in \cite{LafitteMelisSamaey2017} that the spectrum consists of two stationary, well-separated eigenvalue clusters (a fast and slow, dominant cluster). To accommodate these two clusters, the method parameters of projective integration can be selected such that its stability region splits up into two parts.

\begin{enumerate}
\item First, the inner integrator time step $\dt$ is chosen corresponding to the fastest time scale of the problem, which is of the order of $\epsi$. This centers one stability region of the projective method around the fast eigenvalues.

\item Next, the number of inner integrator time steps $K$ is chosen such that all fast eigenvalues lie inside this stability region. In \cite{LafitteMelisSamaey2017}, it was proven that we require $K \ge 2$.

\item Last, the outer integrator time step $\Dt$ is selected such that all dominant eigenvalues fall into the second stability region of the projective method.
\end{enumerate}
Since both $K$ and $\Dt$ are independent of the small-scale parameter $\epsi$, the resulting projective method has a cost that is also independent of $\epsi$, which becomes increasingly advantageous for $\hydroLimit$.

\subsubsection{Continuously spread spectrum}
When considering the linearized BGK equation with $\nu = \rho$, which was studied in \cite{MelisSamaey2017}, the spectrum varies continuously over the negative real axis. This also holds true for the linearized Boltzmann equation, be it on a part of the negative real axis, see figure \ref{fig:spectrumBoltzFourier}. In this case, we require that the stability region of the numerical method does not split up but instead comprises the entire negative real axis up to the fastest eigenvalue of the problem (a numerical method with this property is termed $\mathit{[0,1]}$\textit{-stable}). Here, for simplicity, we assume that the fastest eigenvalue at $t = 0$ corresponds to the fastest possible eigenvalue for all other times $t > 0$. Since [0,1]-stable projective integration methods lose practically all of their potential speed-up, [0,1]-stable telescopic projective integration methods can be designed with much higher speed-ups.

\begin{enumerate}
\item Similarly to projective integration, the innermost integrator time step $h_0$ of the telescopic projective integration method is chosen corresponding to the fastest time scale, which is of the order of $\epsi/\max_{x}\rho(x,0)$.

\item Next, we fix the outermost time step we would like to use, taking into account a CFL-like stability constraint, as follows: $h_L = C\dx$.

\item Before choosing the number of projective levels, we decide on the number of inner integrator time steps $K$, which we consider to be fixed on each projective level. For each chosen value of $K$ there is a corresponding maximal value of $M$ such that the stability region does not split up, see \cite{Gear2003telescopic} or \cite{MelisSamaey2017}.

\item The required number of projective levels $L$ to obtain a [0,1]-stable telescopic method is computed as (see \cite{MelisSamaey2017}):
\begin{equation} \label{eq:tpi_number_of_levels}
	L \approx \frac{\log(h_L) + \log(1/h_0)}{log(M+K+1)}.
\end{equation}

\item For the given values of $h_0$, $h_L$, $K$ and $L$ adapt the value of $M$ on the different projective levels such that the following equation:
\begin{equation} \label{eq:tpi_hL}
	h_L = \prod_{\ell=0}^{L-1} (M_\ell + K + 1)h_0
\end{equation}
is valid (for further details, we refer the reader to \cite{MelisSamaey2017}).
\end{enumerate}
For a [0,1]-stable telescopic method, the values of $h_L$, $M_\ell$ and $K$ are independent of $\epsi$. However, as indicated by equation \eqref{eq:tpi_number_of_levels}, the number of projective levels increases as $O(\log(1/\epsi))$. As a consequence, the cost of a [0,1]-stable telescopic projective integration method is not completely $\epsi$-independent. However, the dependence is rather modest.
  
\subsubsection{Speedup} \label{subsubsec:speedup}
If we assume, as in \cite{Gear2003telescopic}, that the overhead due to extrapolations is negligible and timestepping with the innermost integrator is computationally most demanding, the speedup $\cS_L$ realized by the overall level-$L$ TPI method compared to naive forward Euler timestepping is given by:
\begin{equation} \label{eq:tpi_speedup} 
	\cS_L = \prod_{\ell=0}^{L-1} \frac{M_\ell + K_\ell + 1}{K_\ell + 1},
\end{equation}
that is, the ratio of the total number of naive forward Euler time steps within one outermost time step $h_L$ (see equation \eqref{eq:tpi_hL}) over the number of actual innermost steps in the TPI method.


\section{Numerical experiments} \label{sec:results}
Here, we report simulation results for the BGK equation \eqref{eq:bgk_equation} and the Boltzmann equation \eqref{eqCollision} using pseudo-Maxwellian particles (that is, $\gamma = 0$ and $b_0$ constant). For each experiment, we shall compactly indicate the dimensions in space $(D_x)$ and velocity $(D_v)$ by writing ``$D_x$D/$D_v$D" with $D_x,D_v \in \{1,2\}$. We begin with BGK in 1D/1D (Section~\ref{subsec:bgk_1d1d}), and subsequently consider both the BGK and Boltzmann equation in 1D/2D (Section~\ref{subsec:bgk_boltzmann_1d2d}). Thereafter, we target a shock-bubble interaction problem for the BGK equation in 2D/2D (Section~\ref{subsec:shock_bubble_interaction}), and the more intricated Kelvin-Helmoltz-like instability in the same setting (Section~\ref{subsec:KelvinHelmoltz}). As a last experiment, we deal with the full Boltzmann equation in 2D/2D (Section~\ref{subsec:boltzmann_2d2d}).


\subsection{BGK in 1D/1D} \label{subsec:bgk_1d1d}
As a first experiment, we focus on the nonlinear BGK equation \eqref{eq:bgk_equation} in 1D/1D. We consider a Sod-like test case for $x \in [0,1]$ consisting of an initial centered Riemann problem with the following left and right state values:
\begin{equation} \label{eq:bgk_1d_sod}
	\begin{aligned}
		\begin{pmatrix}	\rho_L \\ \vMacroOneD_L \\ T_L \end{pmatrix} = \begin{pmatrix} 1 \\ 0 \\ 1 \end{pmatrix}, \qquad\qquad 
		\begin{pmatrix}	\rho_R \\ \vMacroOneD_R \\ T_R \end{pmatrix} = \begin{pmatrix} 0.125 \\ 0 \\ 0.25 \end{pmatrix}.
	\end{aligned}
\end{equation}
The initial distribution $\fepsi(x,v,0)$ is then chosen as the Maxwellian \eqref{eq:maxwellian} corresponding to the above initial macroscopic variables. We impose outflow boundary conditions and perform simulations for $t \in [0,0.15]$. 
As velocity space, we take the interval $[-8,8]$, which we discretize on a uniform grid using $J=80$ velocity nodes. In all simulations, space is discretized using the WENO3 spatial discretization with $\dx = 0.01$. Below we regard three gas flow regimes, $\epsi = 10^{-1}$ (kinetic regime), $\epsi = 10^{-2}$ (transitional regime) and $\epsi = 10^{-5}$ (fluid regime), and for each regime, we compare solutions for two cases of collision frequency $\collFreq$ in the BGK equation \eqref{eq:bgk_equation}, $\collFreq = 1$ and $\collFreq = \rho$.

\begin{figure}
	\begin{center}
		\inputtikz{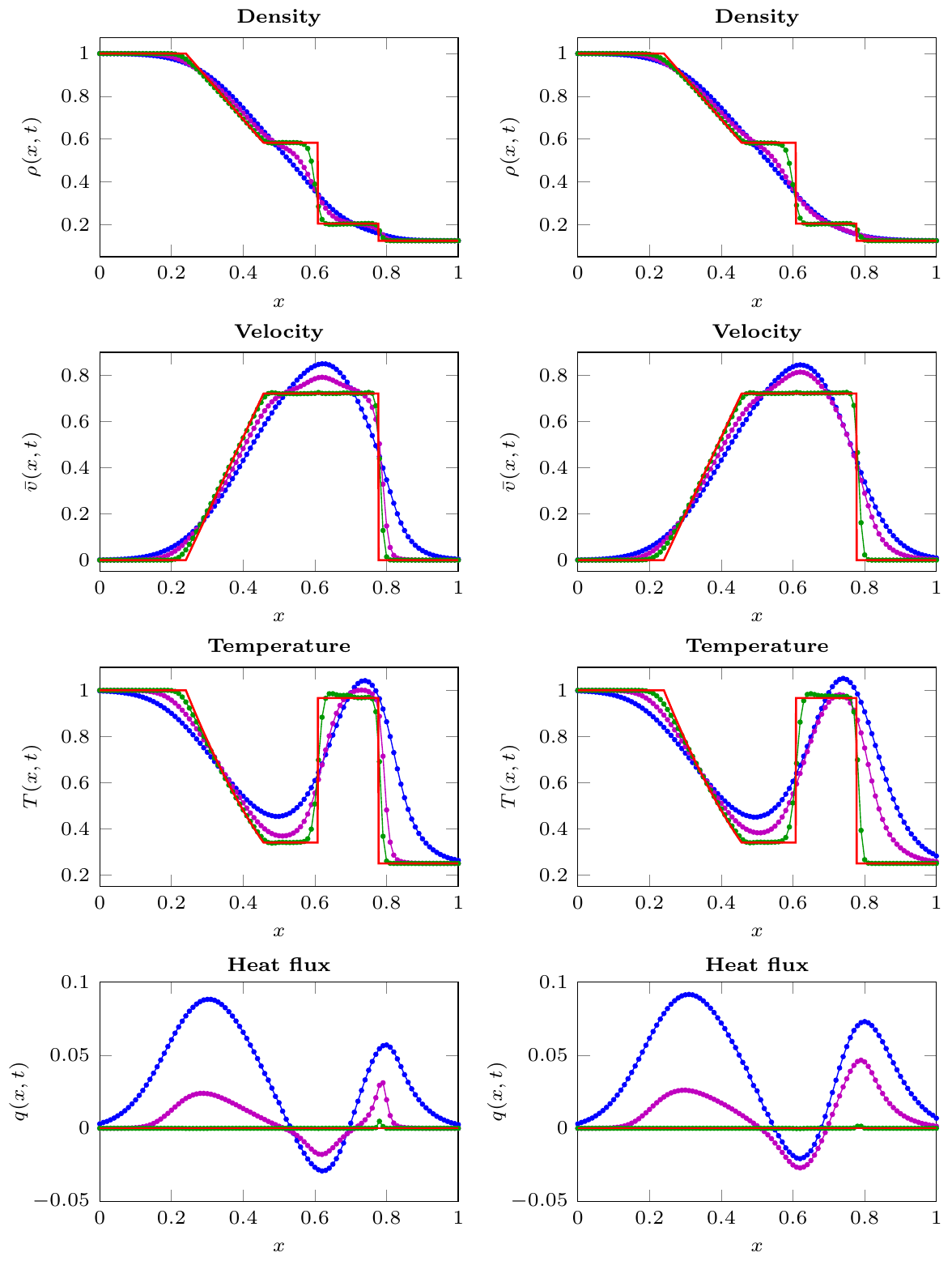}
	\end{center}
	\caption{\label{fig:bgk_1d1d_nu} Numerical solution of the BGK equation in 1D/1D with $\collFreq = 1$ (left) and $\collFreq = \rho$ (right) at $t = 0.15$ for a Sod-like shock test \eqref{eq:bgk_1d_sod} using the WENO3 scheme with $\dx = 0.01$. RK4 is used for $\epsi = 10^{-1}$ (blue dots) and $\epsi = 10^{-2}$ (purple dots) with $\dt = 0.1\dx$. The PRK4 (left) and level-2 TPRK4 (right) methods are used for $\epsi = 10^{-5}$ (green dots). Red line: hydrodynamic limit solution $(\hydroLimit)$. }
\end{figure}

\vspace{-0.5cm}\paragraph{Direct integration $\boldsymbol{(\epsi = 10^{-1}}$ and $\boldsymbol{\epsi = 10^{-2})}$.} In the kinetic $(\epsi = 10^{-1})$ and transitional $(\epsi = 10^{-2})$ regimes, we compute the numerical solution for $\collFreq = 1$ and $\collFreq = \rho$ using the fourth order Runge-Kutta (RK4) time discretization with time step $\dt = 0.1\dx$. The results are shown in figure \ref{fig:bgk_1d1d_nu} for $\collFreq = 1$ (left) and $\collFreq = \rho$ (right), where we display the density $\rho$, macroscopic velocity $\vMacroOneD$ and temperature $T$ as given in \eqref{eq:f_moments} at $t = 0.15$. In addition, we plot the heat flux $q$, which, in a general $D_v$-dimensional setting, is a vector $\heatflux = \sett{q^d}{d}{D_v}$ with components given by:
\begin{equation*}
	q^d = \frac{1}{2}\int_{\R^{D_v}} \abs{\vpec}^2c^d\fepsi d\v,
\end{equation*}
in which $\vpec = \sett{c^d}{d}{D_v} = \v - \vMacroMultiD$ is the peculiar velocity. The different regimes are shown by blue (kinetic) and purple (transitional) dots. The red line in each plot denotes the limiting $(\hydroLimit)$ solution of each macroscopic variable, which all converge to the solution of the Euler system \eqref{eqHydroClosedEuler} with ideal gas law $P = \rho T$ and heat flux $q = 0$.

\vspace{-0.5cm}\paragraph{Projective integration $\boldsymbol{(\epsi = 10^{-5}}$ and $\boldsymbol{\collFreq = 1)}$.} In the fluid regime $(\epsi = 10^{-5})$, direct integration schemes such as RK4 become too expensive due to a severe time step restriction, which is required to ensure stability of the method. Exploiting that the spectrum of the linearized BGK equation with $\collFreq = 1$ resembles that of the linear kinetic models used in \cite{LafitteMelisSamaey2017}, see section \ref{sec:linearOpSpecProp}, we construct a projective integration method to accelerate time integration in the fluid regime. As inner integrator, we select the forward Euler time discretization with $\dt = \epsi$. As outer integrator, we choose the fourth-order projective Runge-Kutta (PRK4) method, using $K=2$ inner steps and an outer step of size $\Dt = 0.4\dx$. Figure \ref{fig:bgk_1d1d_nu} (left) shows the macroscopic observables in the fluid regime for $\collFreq = 1$ at $t = 0.15$ (green dots). From this, we observe that the BGK solution is increasingly dissipative for increasing values of $\epsi$ since the rate with which $\fepsi$ converges to its equilibrium $\MaxBGKOneD$ becomes slower. In contrast, for sufficiently small $\epsi$, relaxation to thermodynamic equilibrium occurs practically instantaneous and the Euler equations \eqref{eqHydroClosedEuler} yield a valid description. Since this is a hyperbolic system, it allows for the development of sharp discontinuous and shock waves which are clearly seen in the numerical solution.

In this numerical test, the speed-up factor between a naive RK4 implementation and the projective integration method is $130.3$, namely formula \eqref{eq:tpi_speedup} with $L=1$ (1 projective level), $K_0=2$ and $M_0=397$.

\begin{figure}
\begin{center}
		\inputtikz{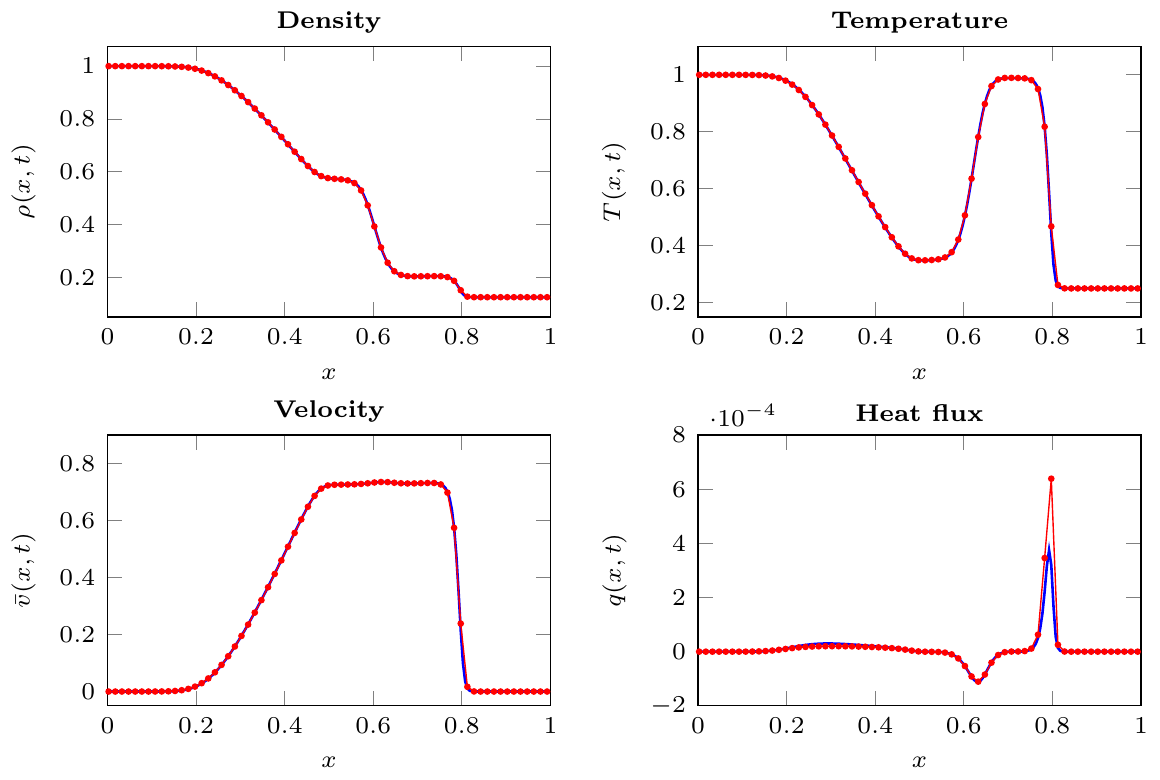}
	\end{center}
	\caption{\label{fig:bgk_rk4TPI_comp} Numerical solution of the BGK equation in 1D/1D with $\collFreq = \rho$, $\epsi = 10^{-5}$   at $t = 0.15$ for a Sod-like shock test \eqref{eq:bgk_1d_sod} using the WENO3 scheme with $\dx = 0.01$. Comparison between level-2 TPRK4 with $\dt = 0.1\dx$ (solid blue line) and classical RK4 with a  $\Dt = 0.5\epsi$ (red dots).}
\end{figure}

\vspace{-0.5cm}\paragraph{Telescopic projective integration $\boldsymbol{(\epsi = 10^{-5}}$ and $\boldsymbol{\collFreq = \rho)}$.} Next, we repeat the above experiment taking $\collFreq = \rho$ in the BGK equation \eqref{eq:bgk_equation}. We now design a $[0,1]$-stable telescopic projective integration method as in \cite{MelisSamaey2017} since, for this choice of $\collFreq$, the spectrum of the linearized BGK equation is spread along the negative real axis and is time-dependent. Therefore, the two-scale nature in case of $\collFreq = 1$ has become a multi-scale problem, which destroys the acceleration in time of projective integration. We construct a $[0,1]$-stable TPRK4 method consisting of $2$ projective levels with FE as innermost integrator with time step $h_0 = \epsi$, constant $K=6$ on each level and an outermost time step $h_2 = 0.4\dx$. The extrapolation step sizes $M$ on each level are calculated as $M=\{14.24, 11.83\}$. The results are shown by green dots in figure \ref{fig:bgk_1d1d_nu} (right). We conclude that the effect of choosing $\collFreq = \rho$ primarily manifests itself in the transitional regime $(\epsi = 10^{-2})$, for which the relaxation rate is not too slow nor too fast. Moreover, it is seen that this choice of collision frequency does not alter the hydrodynamic limit of the BGK equation, which is captured correctly by the telescopic scheme.

Finally, figure \ref{fig:bgk_rk4TPI_comp} compares the solutions obtained with the level-2 TPRK4 method and a classical RK4 method with a very small time step (of order $\epsi$) for the stiff test case where $\epsi =10^{-5}$. We observe a very good agreement between the two simulations (only a very small difference can be seen in the heat flux), while the TPRK4 scheme is more than 10 times faster than the RK4 scheme, because of its bigger time steps. The former scheme would be even more efficient with smaller values of the relaxation parameter.

In this test, the speed-up factor between a naive RK4 implementation and the telescopic projective integration method is $8.2$, namely formula \eqref{eq:tpi_speedup} with $L=2$ (2 projective levels), $(K_0, K_1)=(6,6)$ and $(M_0,M1)=(14.24,11.83)$.

\subsection{BGK and Boltzmann in 1D/2D} \label{subsec:bgk_boltzmann_1d2d}
The BGK equation was introduced as a simplified model for the Boltzmann equation capturing most essential features of the latter. Here, we investigate the difference between both models. Since the Boltzmann collision operator vanishes for a one-dimensional velocity space, in this section, we consider both models in 1D/2D. In the experiments, this is achieved by discretizing space $(x,y)$ on a grid of size $I_x \times 2$ and using homogeneous data along the $y$-direction such that the spatial derivative $\partial_y\fepsi$ exactly cancels out.

We perform the Sod test \eqref{eq:bgk_1d_sod} of the previous section in 1D/2D, see also \cite{Filbet2010}.
As velocity space, we take the domain $[-8,8]^2$, which we discretize on a uniform grid using $J_x = J_y = 32$ velocity nodes along each dimension. In all simulations, space is discretized using the WENO2 spatial discretization with $\dx = 0.01$. Below we regard two regimes, $\epsi = 10^{-2}$ (transitional regime) and $\epsi = 10^{-5}$ (fluid regime), and for each regime, we compare solutions for BGK with $\collFreq = 1$, BGK with $\collFreq = \rho$ and Boltzmann with pseudo-Maxwellian particles. To approximate the Boltzmann collision operator, we apply the fast spectral method described in section \ref{sub:FastSpectral} using $N_\theta = 4$ discrete angles. This is enough because of the spectral accuracy of the trapezoidal rule applied to periodic functions (see \cite{FiMoPa2006} for more details on this topic).

\begin{figure}
	\begin{center}
		\inputtikz{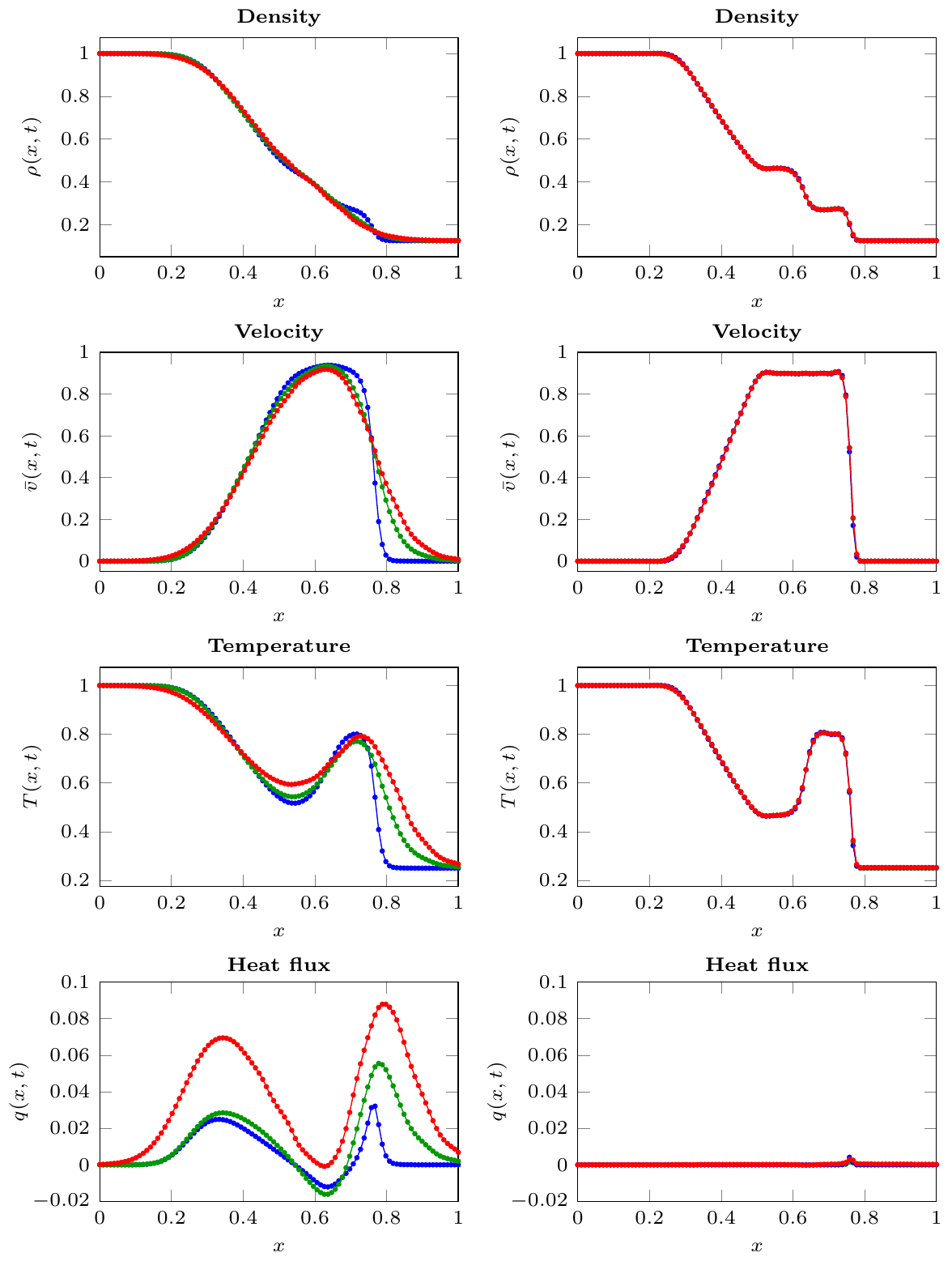}
	\end{center}
	\caption{\label{fig:boltz_bgk_1d2d_nu} Comparison between BGK and Boltzmann in 1D/2D for $\epsi = 10^{-2}$ (left) and $\epsi = 10^{-5}$ (right) at $t = 0.15$ for a Sod-like shock test \eqref{eq:bgk_1d_sod} using the WENO2 scheme with $\dx = 0.01$. Blue dots: BGK with $\collFreq = 1$; green dots: BGK with $\collFreq = \rho$; red dots: Boltzmann. For $\epsi = 10^{-2}$, we used RK4 with $\dt = 0.1\dx$. For $\epsi = 10^{-5}$, we applied a PRK4 (blue dots) and a level-2 TPRK4 (green and red dots) method. }
\end{figure}

\vspace{-0.5cm}\paragraph{Direct integration ($\boldsymbol{\epsi = 10^{-2}}$).} In the transitional regime, we perform all simulations using the RK4 method with time step $\dt = 0.1\dx$, for which we display the results in figure \ref{fig:boltz_bgk_1d2d_nu} (left) for BGK with $\collFreq = 1$ (blue dots), BGK with $\collFreq = \rho$ (green dots) and the Boltzmann equation (red dots). From this, we observe that the BGK solution with $\collFreq = \rho$ is closer to the Boltzmann solution than the BGK solution with $\collFreq = 1$. This is as expected, since the BGK equation with $\collFreq = \rho$ correctly captures the loss term of the Boltzmann collision operator, see \eqref{eq:collision_operator_maxwellian_particles}. Moreover, the discrepancy between Boltzmann and BGK with $\collFreq = \rho$ increases for higher order moments of $\fepsi$; while the density (zeroth order moment) appears to coincide (to the naked eye), the heat flux (third order moment) reveals a clear difference between both models.

\vspace{-0.5cm}\paragraph{Projective methods ($\boldsymbol{\epsi = 10^{-5}}$).} In the fluid regime, the RK4 method becomes too expensive. To that end, for BGK with $\collFreq = 1$, we design a PRK4 method with FE as inner integrator using $\dt = \epsi$, $K = 2$ inner steps and $\Dt = 0.4\dx$. Due to the multi-scale nature of both the BGK relaxation operator with $\collFreq = \rho$ and the Boltzmann collision operator, we construct a $[0,1]$-stable level-2 TPRK4 method for both models using the FE scheme as innermost integrator with $h_0 = \epsi$. We set $K = 4$ constant on each level, compute the extrapolation step sizes as $M = \{14.24, 11.83\}$, and choose the outermost time step as $h_2 = 0.4\dx$. The results can be seen in figure \ref{fig:boltz_bgk_1d2d_nu} (right) using the same plotting style as in the left column. For all models, the projective and telescopic projective integration methods display the expected hydrodynamic limit.

In these numerical tests, the speed-up factor between a naive RK4 implementation and the telescopic projective integration method for the PRK4 method for the BGK model with constant relaxation is $133.3$, namely formula \eqref{eq:tpi_speedup} with $L=1$ (1 projective level), $K_0=2$ and $M_0=397$.
The speedup for the TPRK4 method for the BGK model with nonconstant relaxation rate and the Boltzmann equation is $13$, namely formula \eqref{eq:tpi_speedup} with $L=2$ (2 projective levels), $(K_0, K_1)=(4,4)$ and $(M_0,M1)=(14.24,11.83)$.

\subsection{Shock-bubble interaction in 2D/2D} \label{subsec:shock_bubble_interaction}
\begin{figure}[t]
	\begin{center}
		\inputtikz{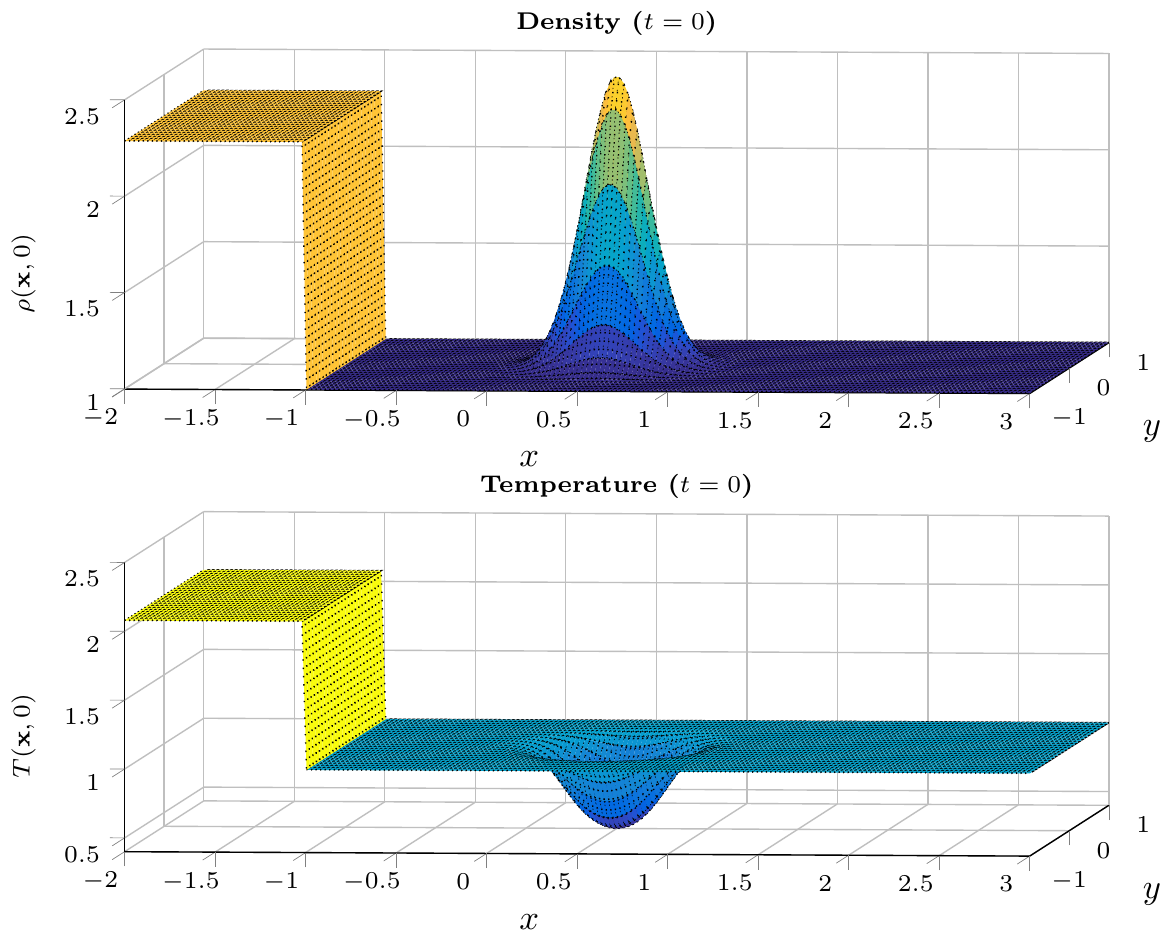}
	\end{center}
	\caption{\label{fig:bgk_2d2d_shock_bubble_t00} Initial solution for density (top) and temperature (bottom) of the shock-bubble interaction in 2D/2D. The spatial domain $[-2,3] \times [-1,1]$ is discretized using $I_x = 200$ and $I_y = 25$. }
\end{figure}

Here, we consider the BGK equation in 2D/2D with constant collision frequency $\collFreq = 1$ and we investigate the interaction between a moving shock wave and a stationary smooth bubble, which was proposed in \cite{Torrilhon2006}, see also \cite{Cai2010}. This problem consists of a shock wave positioned at $x = -1$ in a spatial domain $\x = (x,y) \in [-2,3] \times [-1,1]$ traveling with Mach number $\Ma = 2$ into an equilibrium flow region. Over the shock wave, the following left $(x \le -1)$ and right $(x > -1)$ state values are imposed \cite{Cai2010}:
\begin{equation} \label{eq:bgk_2d_shock}
	\begin{aligned}
		\big(\rho_L, \vMacroOneD^x_L, \vMacroOneD^y_L, T_L\big) &= \left(\frac{16}{7}, \sqrt{\frac{5}{3}}\frac{7}{16}, 0, \frac{133}{64}\right) \\
		\big(\rho_R, \vMacroOneD^x_R, \vMacroOneD^y_R, T_R\big) &= (1, 0, 0, 1).
	\end{aligned}
\end{equation}
Due to this initial profile, the shock wave will propagate rightwards into the flow region at rest $(x > -1)$. Moreover, in this equilibrium region, a smooth Gaussian density bubble centered at $\x_0 = (0.5,0)$ is placed, given by:
\begin{equation} \label{eq:bgk_2d_bubble}
	\rho(\x,0) = 1 + 1.5\exp\left(-16\abs{\x - \x_0}^2\right).
\end{equation}
The initial density $\rho(\x,0)$ and temperature $T(\x,0)$ are visualized in figure \ref{fig:bgk_2d2d_shock_bubble_t00}.
Then, the initial distribution $\fepsi(\x,\v,0)$ is chosen as the Maxwellian \eqref{eq:maxwellian} corresponding to the initial macroscopic variables in \eqref{eq:bgk_2d_shock}-\eqref{eq:bgk_2d_bubble}. We impose outflow and periodic boundary conditions along the $x$- and $y$-directions, respectively, and we perform simulations for $t \in [0,0.8]$.
As velocity space, we take the domain $[-10,10]^2$, which we discretize on a uniform grid using $J_x = J_y = 30$ velocity nodes along each dimension. We discretize space using the WENO2 spatial discretization with $I_x = 200$ and $I_y = 25$. Furthermore, we consider a fluid regime by taking $\epsi = 10^{-5}$.

\begin{figure}
	\begin{center}
		\inputtikz{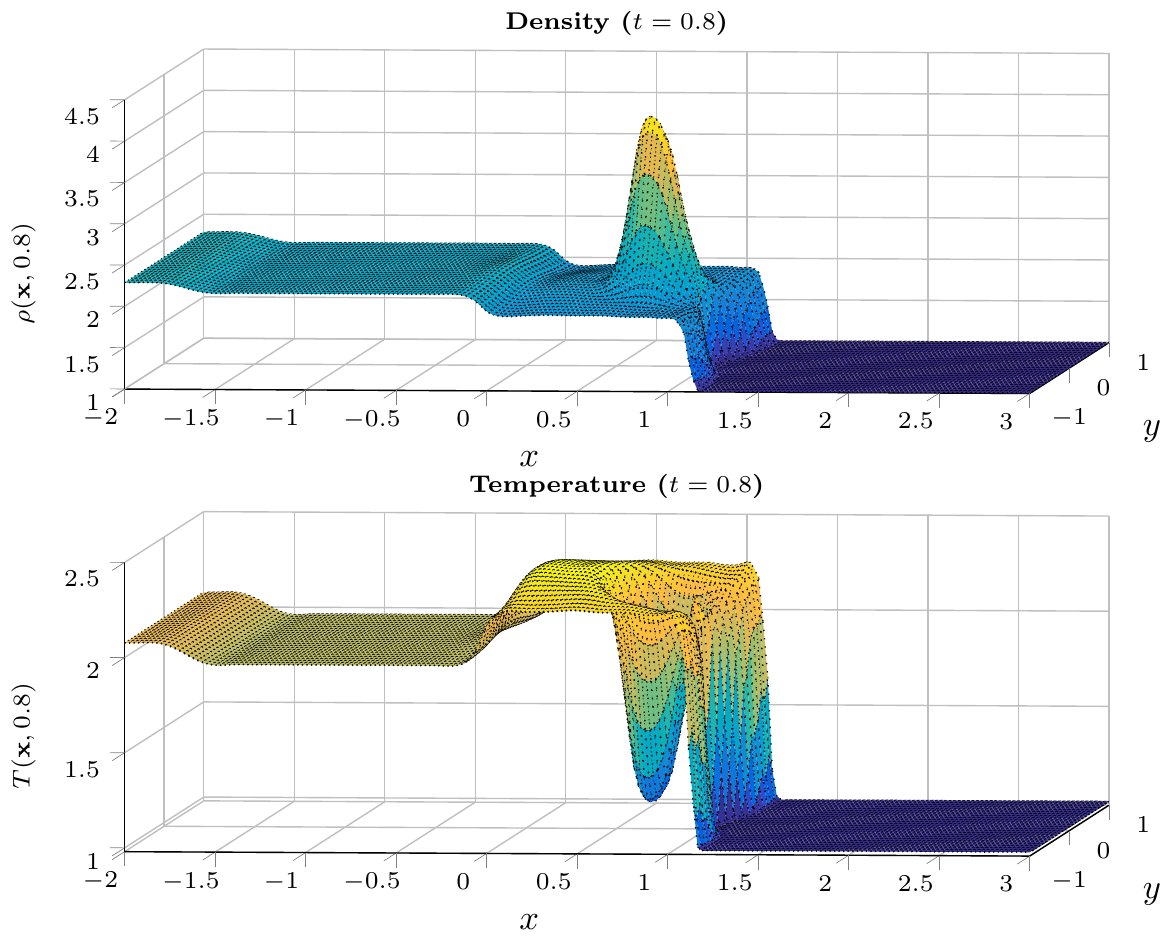}
	\end{center}
	\caption{\label{fig:bgk_2d2d_shock_bubble_t08} Numerical solution of the shock-bubble interaction at $t = 0.8$ using the BGK equation with $\collFreq = 1$ in 2D/2D. We discretized velocity space using $J_x = J_y = 30$. We applied a PRK4 method with FE as inner integrator and $\dt = \epsi = 10^{-5}$ together with the WENO2 spatial discretization scheme with $I_x = 200$ and $I_y = 25$. }
\end{figure}

Since we regard the BGK equation with constant collision frequency $\collFreq = 1$, the spectrum of the linearized BGK operator consists of two eigenvalue clusters. Therefore, we construct a projective integration method to speed up simulation in time. We select the PRK4 method with FE as inner integrator. The inner time step is fixed as $\dt = \epsi$ and we use $K = 2$ inner steps in each outer integrator iteration. The outer time step is chosen as $\Dt = 0.4\dx$. The simulated density and temperature at time $t = 0.8$ are displayed in figure \ref{fig:bgk_2d2d_shock_bubble_t08}. We observe that the shock propagates in the positive $x$-direction and bumps into the stationary bubbly density. A similar behavior is seen for the temperature evolution. 

\begin{figure}
	\begin{center}
		\inputtikz{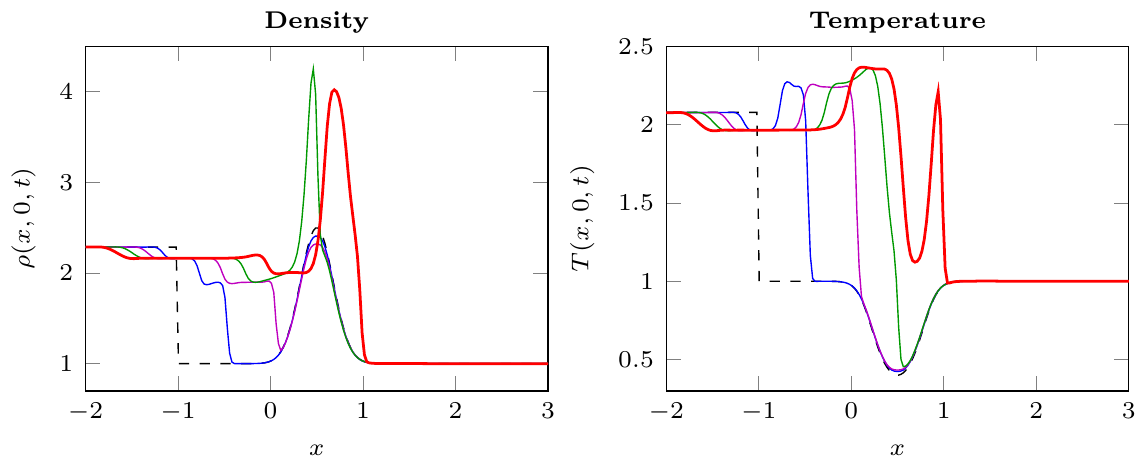}
	\end{center}
	\caption{\label{fig:bgk_2d2d_shock_bubble_y0} Numerical solution of the shock-bubble interaction along $y = 0$ at $t = 0$ (black dashed), $t = 0.2$ (blue), $t = 0.4$ (purple), $t = 0.6$ (green) and $t = 0.8$ (red). }
\end{figure}

To compare our results with those in \cite{Torrilhon2006}, where the smallest value of $\epsi$ is chosen as $\epsi = 10^{-2}$, we regard the one-dimensional evolution of density and temperature along the axis $y = 0$. For $t \in \{0, 0.2, 0.4, 0.6, 0.8\}$, we plot these intersections in figure \ref{fig:bgk_2d2d_shock_bubble_y0}. We conclude that we obtain the same solution structure at $t = 0.8$ as in \cite{Torrilhon2006}. However, our results are sharper and less dissipative supposedly due to the particular small value of $\epsi$ ($10^{-5}$ versus $10^{-2}$). In contrast to \cite{Cai2010}, we nicely capture the swift changes in the temperature profile for $x \in [0.5,1]$ at $t = 0.8$.

Again, in this numerical test, the speed-up factor between a naive RK4 implementation and the projective integration method is $133.3$, namely formula \eqref{eq:tpi_speedup} with $L=1$ (1 projective level), $K_0=2$ and $M_0=397$.

\subsection{A Kelvin-Helmoltz like instability problem} \label{subsec:KelvinHelmoltz}
Keeping the same setting of the BGK equation in 2D/2D with constant collision frequency $\collFreq = 1$, we now turn to a less common test case in the field of collisional kinetic equation, the so-called Kelvin-Helmoltz instability. This phenomenon occurs when two fluids of different densities and in thermodynamic equilibrium move at different speeds. 
It is very well known that such a system will exhibit turbulent, unstable vortices at the interface between the two fluids, because of the velocity shear \cite{Helmoltz1868}. 
In order for these  instabilities to develop, the Reynolds number of the fluids considered must be large. Using the von Karman relation \cite{sone2007molecular}, which states that the Reynolds number is inversely proportional to the Knudsen number $\epsi$, we shall then choose a very small Knudsen number $\epsi = 5\cdot10^{-5}$ along with the following initial condition inspired from \cite{McNallyLyraPassy2012}:
\begin{equation} \label{eq:BGK2dKelvinHelmoltz}
    \begin{pmatrix}	\rho_1 \\ \vMacroOneD^x_1 \\ \vMacroOneD^y_1 \\ T_1 \end{pmatrix} = \begin{pmatrix} 1 \\ 0.5 \\ 0.01 \sin\left (4 \pi x\right ) \\ 1 \end{pmatrix} \quad (y \geq 0), \qquad\quad
	\begin{pmatrix}	\rho_2 \\ \vMacroOneD^x_2 \\ \vMacroOneD^y_2 \\ T_2 \end{pmatrix} = \begin{pmatrix} 2 \\ -0.5 \\ 0.01 \sin\left (4 \pi x\right )  \\ 1 \end{pmatrix} \quad (y < 0).
\end{equation}
The initial distribution $\fepsi(\x,\v,0)$ is chosen as the Maxwellian \eqref{eq:maxwellian} corresponding to the initial macroscopic variables in \eqref{eq:BGK2dKelvinHelmoltz}. We impose periodic and outflow boundary conditions along the $x$- and $y$-directions, respectively, and we perform simulations for $t \in [0,1.6]$.
As velocity space, we take the domain $[-8,8]^2$, which we discretize on a uniform grid using $J_x = J_y = 30$ velocity nodes along each dimension. We discretize space using the WENO2 method on $[-0.5,0.5]\times[-0.5,0.5]$ with $I_x = I_y= 100$.

\begin{figure}
	\begin{center}
		\inputtikz{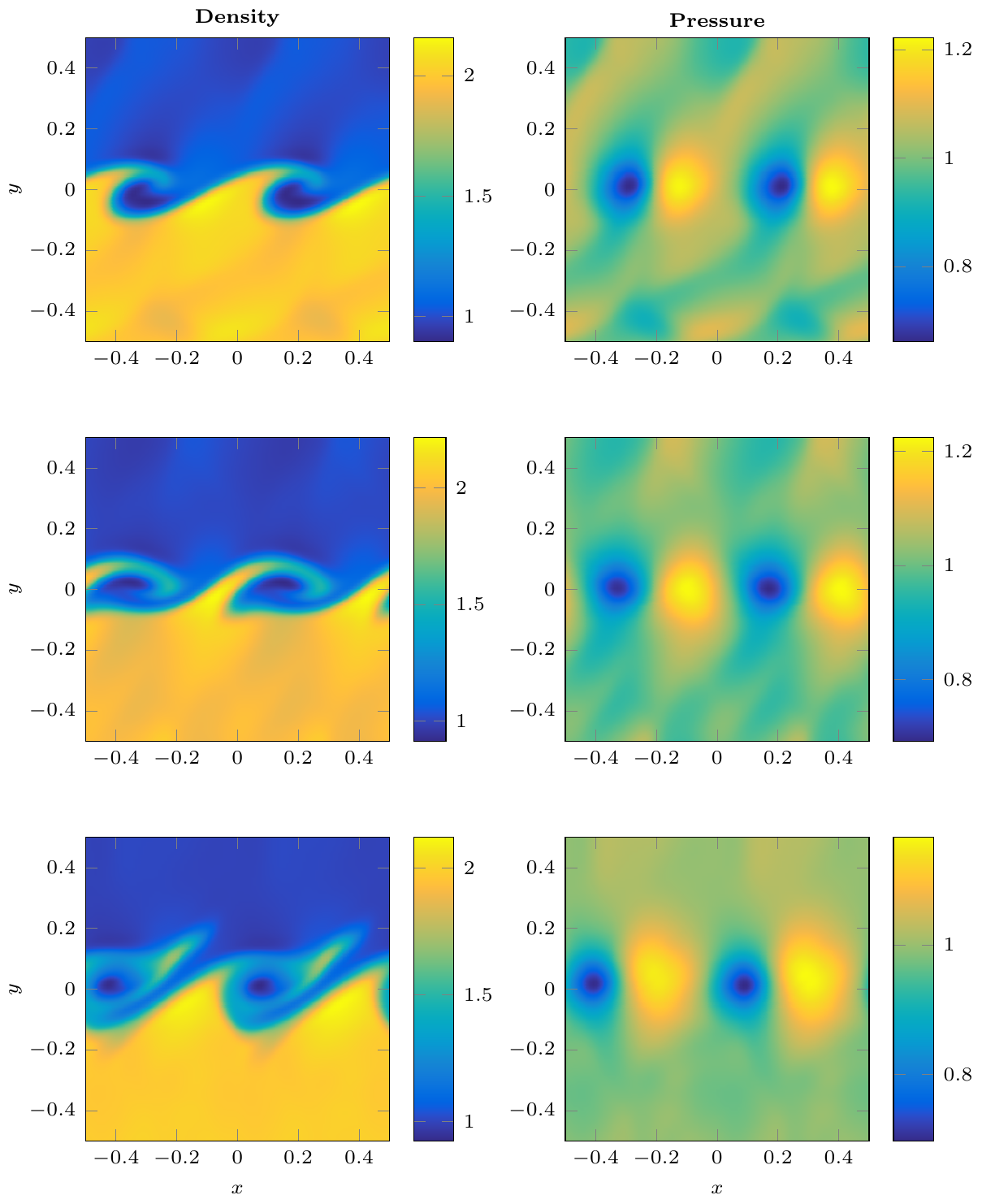}
	\end{center}
	\caption{\label{fig:bgk_2d2d_KH_Density-Pressure} Density (left) and pressure (right) of the Kelvin-Helmoltzm-like instability at times $t=0.6$ (first row), $0.9$ (second row) and $1.6$ (third row). }
\end{figure}

\begin{figure}
	\begin{center}
		\inputtikz{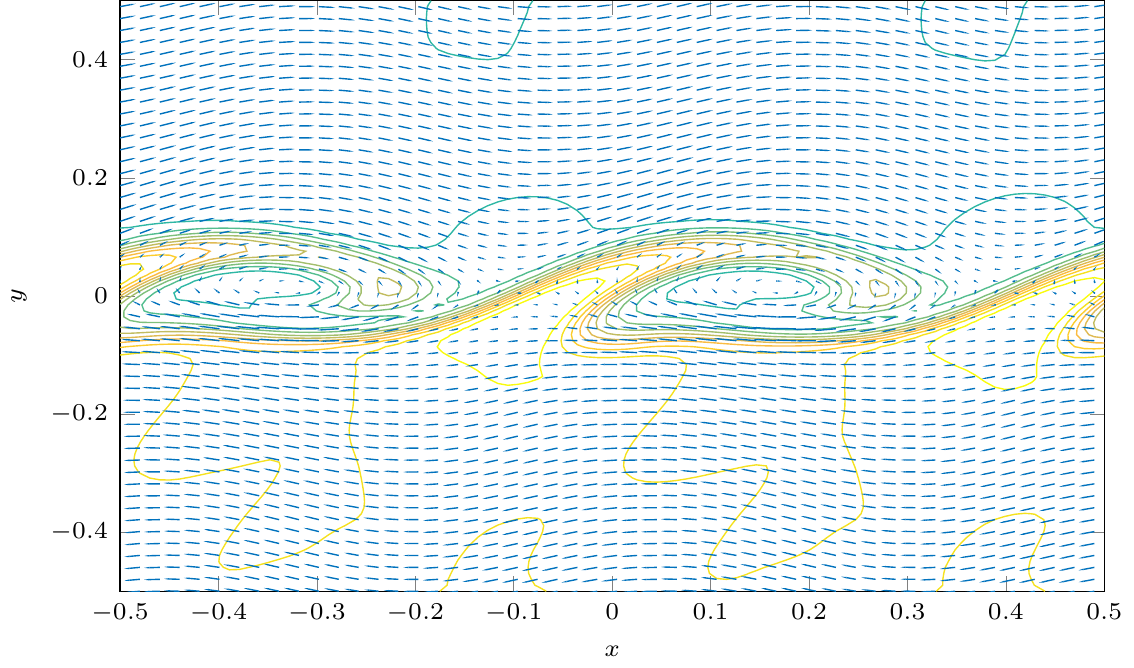}
	\end{center}
	\caption{\label{fig:bgk_2d2d_KH_velocityField}Velocity field and density line of the Kelvin-Helmoltzm-like instability at time $t=0.9$.}
\end{figure}

Since we consider again the BGK equation with constant collision frequency $\collFreq = 1$, the spectrum of the linearized BGK operator consists of two eigenvalue clusters. Therefore, we construct a projective integration method to speed up simulation in time. We select the PRK4 method with FE as inner integrator. The inner time step is fixed as $\dt = \epsi$ and we use $K = 3$ inner steps in each outer integrator iteration (this test case is stiffer than the previous one because of the turbulent regime). The outer time step is chosen as $\Dt = 0.45\dx$. 
The simulated density and pressure at time $t = 0.4$, $0.9$ and $1.6$ are displayed in figure \ref{fig:bgk_2d2d_KH_Density-Pressure}, along with the vector field $ \vMacroMultiD$ with the contour lines of density at time $t=0.9$ in figure \ref{fig:bgk_2d2d_KH_velocityField}.
  
In this numerical test, the speed-up factor between a naive RK4 implementation and the projective integration method is $22.5$, namely formula \eqref{eq:tpi_speedup} with $L=1$ (1 projective level), $K_0=3$ and $M_0=86$.

As expected, we observe vortices developing along the velocity shear line $y=0$, and then expanding, forming small scale structures as time evolves, hence validating the turbulent behavior exhibited by the fluid for this Reynolds number. Note that the use of high-order methods in both space (WENO2) and time (PRK4) is crucial, because of these fine scale structures. Choosing less accurate (and more dissipative) schemes, such as forward Euler in time and upwind in space, leads to a uniform density and pressure instead of vortices in our numerical tests. 

\subsection{Boltzmann in 2D/2D} \label{subsec:boltzmann_2d2d}
As a last experiment, we concentrate on the Boltzmann equation with pseudo-Maxwellian particles in 2D/2D. As initial configuration for the gas, we consider the double Sod shock test, that is, for ${\x = (x,y) \in [-0.5,0.5]^2}$, we set:
\begin{equation} \label{eq:boltzmann_2d_sod}
	\begin{aligned}
		\begin{pmatrix}	\rho_1 \\ \vMacroMultiD_1 \\ T_1 \end{pmatrix} = \begin{pmatrix} 0.1 \\ \boldsymbol{0} \\ 1 \end{pmatrix} \quad (xy \le 0), \qquad\quad
		\begin{pmatrix}	\rho_2 \\ \vMacroMultiD_2 \\ T_2 \end{pmatrix} = \begin{pmatrix} 1 \\ \boldsymbol{0} \\ 1 \end{pmatrix} \quad (\text{otherwise}).
	\end{aligned}
\end{equation}
The initial distribution $\fepsi(\x,\v,0)$ is then chosen as the Maxwellian \eqref{eq:maxwellian} corresponding to the above macroscopic variables. We impose outflow boundary conditions along both dimensions and perform simulations for $t \in [0,0.16]$. 
As velocity space, we take the domain $[-8,8]^2$, which we discretize on a uniform grid with $J_x = J_y = 32$ velocity nodes along each dimension. Furthermore, we discretize the spatial domain using the WENO2 spatial discretization with $I_x = I_y = 64$ grid points along each dimension, and we fix $\epsi = 5 \cdot 10^{-5}$ (that is, we are in the fluid regime). The Boltzmann collision operator is approximated again using the fast spectral method described in section \ref{sub:FastSpectral} with $N_\theta = 4$ discrete angles.

\begin{figure}
	\begin{center}
		\inputtikz{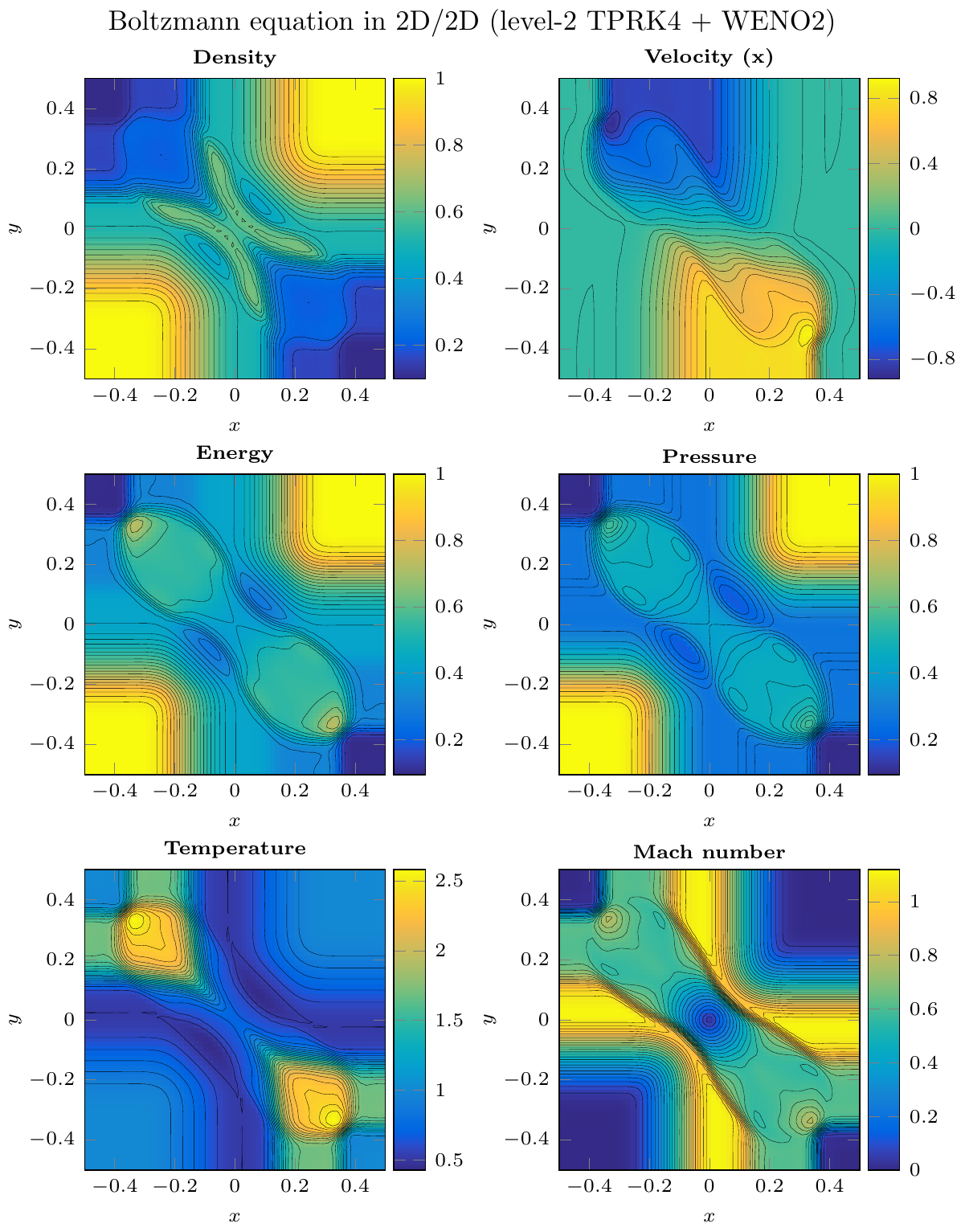}
	\end{center}
	\caption{\label{fig:boltz_2d2d_sod} Numerical solution of the Boltzmann equation with pseudo-Maxwellian particles in 2D/2D at $t = 0.16$ for a double Sod shock test \eqref{eq:boltzmann_2d_sod}. Velocity space is discretized using $J_x = J_y = 32$. We applied a level-2 TPRK4 method with FE as innermost integrator and $h_0 = \epsi = 5 \cdot 10^{-5}$ together with WENO2 with $I_x = I_y = 64$. }
\end{figure}

For the time simulation of the Boltzmann equation, we apply a level-2 TPRK4 method with FE as innermost integrator using $h_0 = \epsi$ as innermost time step. We set $K = 3$ constant on each level and compute the extrapolation step sizes as $M = \{6.66, 4.80\}$. The outermost time step is chosen as $\Dt = 0.3\dx$. In figure \ref{fig:boltz_2d2d_sod}, we plot various macroscopic observables of interest at $t = 0.16$. The density, macroscopic velocity along $x$ and temperature are computed as the moments of $\fepsi$, see \eqref{eq:f_moments}. Then, pressure, energy and the Mach number are obtained, respectively, as:
\begin{equation*}
	P = \rho T, \qquad\quad E = \frac{1}{2}\rho\abs{\vMacroMultiD}^2 + P, \qquad\quad \Ma = \frac{\abs{\vMacroMultiD}}{\sqrt{T}}.
\end{equation*}

In this test, the speed-up factor between a naive RK4 implementation and the telescopic projective integration method is $5.9$, namely formula \eqref{eq:tpi_speedup} with $L=2$ (2 projective levels), $(K_0, K_1)=(3,3)$ and $(M_0,M1)=(6.66,4.80=)$.

\section{Conclusions} \label{sec:conclusions}

We extended projective and telescopic projective integration methods, studied in \cite{LafitteMelisSamaey2017,MelisSamaey2017} for kinetic BGK-like equations, to allow for a fully explicit, high-order time simulation of both the nonlinear BGK equation and the full Boltzmann equation for pseudo-Maxwellian particles. We developed a framework of linearized operators, which revealed that the linearized BGK operator closely resembles the relaxation operator of the linear kinetic models used in \cite{LafitteMelisSamaey2017,MelisSamaey2017}. As a result, the design of stable projective and telescopic projective integration methods for simulation of the nonlinear BGK equation is practically identical to that presented in \cite{LafitteMelisSamaey2017} and \cite{MelisSamaey2017}, respectively. Since the spectrum of the linearized Boltzmann equation for pseudo-Maxwellian particles is known to be spread along the negative real axis, accelerated explicit time integration of the full Boltzmann equation required the construbtion of $[0,1]$-stable TPI methods. Although such methods are not completely asymptotic-preserving, the cost scales only logarithmically with the stiffness of the problem, which we consider acceptable. We designed and applied the projective methods to a number of model problems of increasing dimension and complexity, showing the potential of these schemes, with speed-ups that ranges from 8 for the most complicated Boltzmann model to more than a hundred for the BGK equation with linear relaxation rate.

Because of their explicit nature and hierarchical structure, telescopic projective integration methods naturally lend themselves to adaptivity. In future work, we will explore the use of space-dependent hierarchies of projective levels, adapted to the local values of $\epsi$, to further increase computational efficiency. 

\bibliographystyle{acm}
\bibliography{refs}

\appendix

\section{Free transport with WENO scheme} \label{sub:weno}
In weighted essentially non-oscillatory (WENO) methods, introduced in \cite{Liu1994}, the stencil is built adaptively for each finite volume cell $\cell_i = [x_{i-1/2},x_{i+1/2}]$. A WENO method considers all $k$ possible stencils for cell $\cell_i$ at once and computes a weighted average over the $k$ resulting reconstructed solutions. The $k$ stencils for cell $\cell_i$ are given by:
\begin{equation} \label{eq:weno_stencils}
	\sten_i^\ell = \{\cell_{i-\ell}, \ldots, \cell_{i-\ell+k-1}\} \qquad (\ell = 0, \ldots, k-1),
\end{equation}
in which the left shift $\ell$ is used to iterate over the possible stencils.
As in \cite{Shu1998}, for each stencil in \eqref{eq:weno_stencils}, the reconstructed solution on the left side of $x_{i+1/2}$ and on the right side of $x_{i-1/2}$ are calculated as:
\begin{equation} \label{eq:weno_reconstructed_solutions}
	\begin{aligned}
		u_{i+1/2}^\ell &= \sum_{r=0}^{k-1} c_{\ell,r}U_{i-\ell+r} \\
		u_{i-1/2}^\ell &= \sum_{r=0}^{k-1} c_{\ell-1,r}U_{i-\ell+r}
	\end{aligned}
	\qquad\quad (\ell = 0, \ldots, k-1),
\end{equation}
where we dropped the minus and plus superscripts in \eqref{eq:weno_reconstructed_solutions}, respectively, for clarity. In what follows, we concentrate on the right boundary of $\cell_i$. The left boundary is treated analogously.

Using the $k$ reconstructed solutions in \eqref{eq:weno_reconstructed_solutions}, we could define the following weighted average:
\begin{equation} \label{eq:weno_reconstruction_averaged_smooth}
	u_{i+1/2}^- = \sum_{\ell=0}^{k-1} d_{\ell}u_{i+1/2}^\ell,
\end{equation}
for some weights $d_{\ell}$, $\ell=0, \ldots, k-1$. Equation \eqref{eq:weno_reconstruction_averaged_smooth} can be seen as the reconstruction of the solution on the left side of $x_{i+1/2}$ using a stencil of $2k-1$ consecutive cells. The weights $d_{\ell}$ can be chosen such that \eqref{eq:weno_reconstruction_averaged_smooth} has an order of accuracy of $2k-1$, that is:
\begin{equation*}
	u_{i+1/2}^- = u(x_{i+1/2}) + O\left(\dx^{2k-1}\right),
\end{equation*}
which holds if the solution is smooth over all cells of \eqref{eq:weno_stencils}. In \cite{Arandiga2011} it is proven that the weights $d_{\ell}$ can be written in terms of binomial coefficients as:
\begin{equation*}
	d_{\ell} = \dfrac{\dbinom{k-1}{k-1-\ell}\dbinom{k}{k-1-\ell}}{\dbinom{2k-1}{k}} \qquad (\ell = 0, \ldots, k-1),
\end{equation*}
from which it is straightforward to see that $0 < d_{\ell} \le 1$ and $\sum_{\ell=0}^{k-1} d_{\ell} = 1$ such that \eqref{eq:weno_reconstruction_averaged_smooth} represents a convex combination of the $k$ reconstructed solutions. 
For instance, we compute:
\begin{equation} \label{eq:weno_dl}
	\begin{alignedat}{4}
		k &=1: &&\quad d_0 = 1&& \\
		k &=2: &&\quad d_0 = \frac{2}{3}, &&\quad d_1 = \frac{1}{3} \\
		k &=3: &&\quad d_0 = \frac{3}{10}, &&\quad d_1 = \frac{6}{10}, \quad d_2 = \frac{1}{10}.
	\end{alignedat}
\end{equation}
The idea of WENO methods is to use the weights $d_{\ell}$ only for stencils over which the solution is smooth and allocate very small weights to stencils containing one or more discontinuities. If we denote these WENO weights by $\omega_{\ell}$, we write the WENO reconstruction as:
\begin{equation*} 
	u_{i+1/2}^- = \sum_{\ell=0}^{k-1} \omega_{\ell}u_{i+1/2}^\ell,
\end{equation*}
for which we require that the weights $\omega_{\ell}$ also form a convex set, that is: 
\begin{equation} \label{eq:weno_omegal_convex}
	0 < \omega_\ell \le 1, \qquad\quad
	\sum_{\ell=1}^{k-1}\omega_\ell = 1.
\end{equation}
In \cite{Liu1994} the second constraint in \eqref{eq:weno_omegal_convex} is guaranteed by writing $\omega_\ell$ in terms of constants $\alpha_\ell$ as:
\begin{equation*} 
	\omega_\ell = \frac{\alpha_\ell}{\displaystyle\sum_{s=0}^{k-1} \alpha_s} \qquad (\ell = 0, \ldots, k-1).
\end{equation*}
The constants $\alpha_\ell$ are expressed in terms of the weights $d_{\ell}$, as given in \eqref{eq:weno_dl}, and new constants $\beta_{\ell}$, leading to:
\begin{equation} \label{eq:weno_alphal}
	\alpha_\ell = \frac{d_\ell}{(\delta + \beta_\ell)^2} \qquad (\ell = 0, \ldots, k-1),
\end{equation}
in which the coefficients $\beta_\ell$ determine the smoothness of each stencil and are therefore termed the \emph{smoothness indicators}. Since the smoothness indicators can become zero, a small constant $\delta > 0$ is added in \eqref{eq:weno_alphal} that avoids division by zero. A typical choice is $\delta = 10^{-6}$ and the numerical experiments in \cite{Jiang1996} suggest that the resulting WENO method is not sensitive to the choice of $\delta$ as soon as $10^{-7} \le \delta \le 10^{-5}$. 

The critical step in the design of WENO methods appears to be the definition of the smoothness indicators $\beta_{\ell}$. In \cite{Jiang1996}, the authors proposed the following smoothness indicators:
\begin{equation} \label{eq:weno_smoothness_indicators}
	\beta_\ell = \sum_{s=1}^{k-1}\dx^{2s-1} \int_{x_{i-1/2}}^{x_{i+1/2}} \left(\frac{d^s p_i^\ell(x)}{dx^s}\right)^2 dx \qquad (\ell = 0, \ldots, k-1).
\end{equation}
Each smoothness indicator in \eqref{eq:weno_smoothness_indicators} includes the information of all derivatives that are approximated by the reconstruction polynomial $p_i^\ell(x)$ based on stencil $\sten_i^\ell$ as a measure of smoothness of the true solution in cell $\cell_i$. The factor $\dx^{2s-1}$ in front of the integral ensures that $\beta_\ell$ is independent of $\dx$ when computing the derivative of $p_i^\ell(x)$ in \eqref{eq:weno_smoothness_indicators}. For instance, for $k = 2$, the smoothness indicators in \eqref{eq:weno_smoothness_indicators} are calculated as:
\begin{equation} \label{eq:weno_betal_k2}
	\begin{aligned}
		\beta_0 &= (U_{i} - U_{i+1})^2 \\
		\beta_1 &= (U_{i-1} - U_{i})^2.
	\end{aligned}
\end{equation}
For higher-order reconstructions, the derivatives in \eqref{eq:weno_smoothness_indicators} can be computed by using a general expression of the reconstruction polynomial. For $k=3$, after rearranging terms, we find (see \cite{Jiang1996,Shu1998}):
\begin{equation} \label{eq:weno_betal_k3}
	\begin{aligned}
		\beta_0 &= \frac{1}{4}(3U_{i} - 4U_{i+1} + U_{i+2})^2 + \frac{13}{12}(U_{i} - 2U_{i+1} + U_{i+2})^2 \\
		\beta_1 &= \frac{1}{4}(U_{i-1} - U_{i+1})^2 + \frac{13}{12}(U_{i-1} - 2U_{i} + U_{i+1})^2 \\
		\beta_2 &= \frac{1}{4}(U_{i-2} - 4U_{i-1} + 3U_{i})^2 + \frac{13}{12}(U_{i-2} - 2U_{i-1} + U_{i})^2,
	\end{aligned}
\end{equation}
see also the work of \cite{Henrick2005}. For orders $k=4,5,6$, the corresponding expressions were obtained in \cite{Balsara2000} and extended to even higher orders $k=7,8,9$ in \cite{Gerolymos2009}. From the expressions in \eqref{eq:weno_betal_k2} and \eqref{eq:weno_betal_k3}, we see that $\beta_\ell = O(\dx^2)$ as soon as the solution is smooth over the stencil, while we have $\beta_\ell = O(1)$ for stencils containing a discontinuity. Using equation \eqref{eq:weno_alphal}, this results in weights $\omega_\ell \approx d_\ell$ for smooth solutions and $\omega_\ell = O(\dx^4)$ for discontinuous solutions, as desired.
  
\section{Evaluating the Boltzmann collision operator with a fast spectral scheme} \label{sub:FastSpectral}

The fast spectral discretization of the Boltzmann operator taken from \cite{MoPa:2006} employed in this work is described in this appendix.
 
To this aim, and since the Boltzmann collision operator acts only on the velocity variables, we focus on a given spatial cell $\x_j$ at a given instant of time $t^n$. 
Hence, only the dependency on the velocity variable $\v$ is considered for the distribution function $f$, i.e. $f=f(\v)$.

The first step to construct our spectral discretization is to truncate the integration domain of the Boltzmann integral \eqref{eq:collision_operator}. 
As a consequence, we suppose the distribution function $f$ to have compact support on the ball $\Ball_0(R)$ of
radius $R$ centered in the origin. 
Since one can prove (see e.g. \cite{PaRuSINUM2000}) that $\supp (Q(f)(v)) \subset \Ball_0({\sqrt 2}R)$, in order to write a spectral approximation which avoids aliasing, it is sufficient that the distribution function $f(\v)$ is restricted on the cube $[-T,T]^{D_v}$ with $T \geq (2+{\sqrt 2})R$. Successively, one should assume $f(\v)=0$ on $[-T,T]^{D_v} \setminus \Ball_0(R)$ and extend $f$ to a periodic function on the set $[-T,T]^{D_v}$. 
Let observe that the lower bound for $T$ can be improved. For instance, the choice $T=(3+{\sqrt 2})R/2$ guarantees the absence of intersection between periods where $f$ is different from zero. However, since in practice the support of $f$ increases with time, we can just minimize the errors due to aliasing \cite{canuto:88} with spectral accuracy. 

To further simplify the notation, let us take $T=\pi$ and hence $R=\lambda\pi$ with $\lambda = 2/(3+\sqrt{2})$ in the following. 
We denote by $\QL_B(f)$ the Boltzmann operator with cut-off. Hereafter, using one (bold) index to denote the $D_v$-dimensional sums, we have that the approximate function $f_N$ can be represented as the truncated Fourier
series by
\begin{equation}
f_N(\v) = \sum_{\k=-N/2}^{N/2} \f_\k e^{i \k \cdot \v},
\label{eq:FU}
\end{equation}
where the $\k^{th}$ Fourier coefficient is given by 
\begin{equation*}
\f_\k = \frac{1}{(2\pi)^{D_v}}\int_{[-\pi,\pi]^{D_v}} f(\v)
e^{-i \k \cdot \v }\,d\v.
\end{equation*}
We then obtain a spectral quadrature of our collision operator by projecting \eqref{eq:collision_operator} on the space
of trigonometric polynomials of degree less or equal to $N$, i.e.
\begin{equation}
{\hat \Q}_\k=\int_{[-\pi,\pi]^{D_v}}
\QL_B(f_N) \,
e^{-i \k \cdot \v}\,dv, \quad \k=-N/2,\ldots,N/2. 
\label{eq:VAR}
\end{equation}
Finally, by substituting expression \eqref{eq:FU} in \eqref{eq:VAR} one gets after some computations
\begin{equation}
{\hat \Q}_\k = \sum_{\substack{\l,\m=-N/2\\ \l+\m=\k}}^{N/2} \f_\l\,\f_\m \,
\bb(\l,\m),\quad \k=-N,\ldots,N,
\label{eq:CF1}
\end{equation}
where $\bb(\l,\m)=\B(\l,\m)-\B(\m,\m)$ are given by
\begin{equation*}
\B(\l,\m) = \int_{\Ball_0(2\lambda\pi)}\int_{\mathbb{S}^{D_v-1}} 
B(|q|, \cos\theta) e^{-i(\l\cdot q^++\m\cdot q^-)}\,d\omega\,dq.
\end{equation*}
with \begin{equation*}
q^{+} = \frac12(q+\vert q\vert \omega), \quad
q^{-} = \frac12(q-\vert q\vert \omega).
\end{equation*}
Let us notice that the naive evaluation of (\ref{eq:CF1}) requires $O(n^2)$ operations, where $n=N^3$. This causes the spectral method to be computationally very expensive,
especially in dimension three. In order to reduce the number of operations needed to evaluate the collision integral, the main idea is to use another representation of \eqref{eq:collision_operator}, the so-called Carleman representation \cite{Carl:EB:32} which is obtained by using the  following identity
  \begin{equation*}
    \frac{1}{2} \, \int_{\mathbb{S}^{D_v-1}} F(|u|\sigma - u) \, d\sigma	 = \frac{1}{|u|^{d-2}} \, \int_{\mathbb{R}^{D_v}} \delta(2 \, x \cdot u + |x|^2) \, F(x) \, dx.
 \end{equation*} 
This gives in our context for the Boltzmann integral
\begin{equation*}
 \Q (f)= \int_{\R^{D_v}} \int_{\R^{D_v}} {\tilde B}(x,y) 
 \delta(x \cdot y) 
 \left[ f(v + y) \, f(v+ x) - f(v+x+y) \, f(v) \right] \, dx \,
 dy,
 \end{equation*} 
with 
  \begin{equation}\label{eq:Btilde}
  \tilde{B}(|x|,|y|) =
  2^{D_v-1} \, \sigma\left(\sqrt{|x|^2+|y|^2}, \frac{|x|}{\sqrt{|x|^2+|y|^2}} \right) \, (|x|^2+|y|^2)^{-\frac{D_v-2}2}.
  \end{equation}
This transformation yields the following new spectral quadrature formula 
 \begin{equation}\label{eq:ode}
 \hat{\Q}_\k  =
 \sum_{\underset{\l+\m=\k}{\l,\m=-N/2}}^{N/2} {\hat{\beta}}_F(\l,\m) \, \hat{f}_\l \, \hat{f}_\m, \ \ \
 \k=-N,...,N
 \end{equation}
where ${\hat{\beta}}_F(\l,\m)=\B_F(\l,\m)-\B_F(\m,\m)$ are now 
given by
 \begin{equation*}
 \B_F(\l,\m) = \int_{\Ball_0(R)} \int_{\Ball_0(R)}
 \tilde{B}(x,y) \, \delta(x \cdot y) \, 
 e^{i (\l \cdot x+ \m \cdot y)} \, dx \, dy.
 \end{equation*}
Now, in order to reduce the number of operation needed to evaluate~\eqref{eq:ode}, we look for a convolution structure. 
The aim is to approximate each
${\hat{\beta}}_F(\l,\m)$ by a sum
 \[ {\hat{\beta}}_F(\l,\m) \simeq \sum_{p=1} ^{A} \alpha_p (\l) \alpha' _p (\m), \]
where $A$ represents the number of finite possible directions of collisions.
This finally gives a sum of $A$ discrete convolutions and, consequently, the algorithm can be computed in $O(A \, N \log_2 N)$ operations by means of
standard FFT technique~\cite{canuto:88}. 

In order to get this convolution form, we make the decoupling assumption
 \begin{equation*}
 \tilde{B}(x,y) = a(|x|) \, b(|y|).
 \end{equation*}
This assumption is satisfied if $\tilde B$ is constant. This is the case of Maxwellian molecules in dimension two, which is the case we shall consider for the numerical simulations of section \ref{sec:results}.
Indeed, using the kernel \eqref{eq:collision_operator_maxwellian_particles} in \eqref{eq:Btilde}, one has
\[
  \tilde{B}(x,y) = 2^{D_v - 1} b_0 (|x|^2+|y|^2)^{-\frac{D_v-\alpha-2}2},
\] 
so that $\tilde{B}$ is constant if $D_v = 2$ and $\alpha = 0$.
Here we write $x$ and $y$ in spherical coordinates $x = \rho e$ and $y = \rho' e'$ to get
 \begin{equation*}
 \B_F(\l,\m) = \frac14 \, \int_{\mathbb{S}^1} \int_{\mathbb{S}^1}
 \delta(e \cdot e') \,
 \left[ \int_{-R} ^R e^{i \rho (\l \cdot e)} \, d\rho \right] \,
 \left[ \int_{-R} ^R e^{i \rho' (\m \cdot e')} \, d\rho' \right] \, de \,
 de'.
 \end{equation*}
Then, denoting $ \phi_R ^2 (s) = \int_{-R} ^R e^{i \rho s} \, d\rho,$ for $s \in \R$,  
we have the explicit formula
 \[ \phi_R ^2 (s) = 2 \, R \,{\sinc} (R s), \]
 where ${\sinc}(x)=\frac{\sin(x)}{x}$.
This explicit formula is further plugged in the expression of $\B_F(l,m)$ and using its parity property, this yields
 \begin{equation*}
 \B_F (\l,\m) =  \int_0 ^{\pi} \phi_R ^2 (\l \cdot e_{\theta})\, \phi_R ^2 (\m \cdot e_{\theta+\pi/2}) \,
 d\theta.
 \end{equation*}
Finally, a regular discretization of $N_\theta$ equally spaced points, which is spectrally accurate because of the periodicity of the function \cite{kurganov2007spectral}, gives
 \begin{equation*}
 \B_F (\l,\m) = \frac{\pi}{M} \, \sum_{p=1} ^{N_\theta} \alpha_p (\l) \alpha' _p (\m),
 \end{equation*}
with
 \begin{equation*} 
 \alpha _p (\l) = \phi_R ^2 (\l \cdot e_{\theta_p}), \hspace{0.8cm} \alpha' _p (\m) = \phi_R ^2 (\m \cdot e_{\theta_p+\pi/2}) 
 \end{equation*}
where $\theta_p = \pi p/N_\theta$.

\end{document}